\makeatletter
\def\input@path{{./INFORMS-IJOC-Template-6-10-2024}}
\makeatother

\documentclass[12pt]{article}

\usepackage{ifthen}
\newboolean{ijoc}
\usepackage{natbib}
 \bibpunct[, ]{(}{)}{,}{a}{}{,}%

\ifthenelse{\boolean{ijoc}}{

\usepackage{eqndefns-left} % For checking the display equation width and equation environment definitions %
\RequirePackage{tgtermes}
\RequirePackage{newtxtext}
\RequirePackage{newtxmath}
\RequirePackage{bm}
\RequirePackage{endnotes}

\OneAndAHalfSpacedXII % Current default line spacing
\usepackage{algorithm}
\usepackage{algpseudocode}
\usepackage{tikz}
\EquationsNumberedThrough    % Default: (1), (2), ...
\TheoremsNumberedThrough     % Preferred (Theorem 1, Lemma 1, Theorem 2)
\ECRepeatTheorems  %  

\MANUSCRIPTNO{}

\RUNAUTHOR{Zhu and Burer}
\RUNTITLE{An Extended Validity Domain for Constraint Learning}
\TITLE{An Extended Validity Domain for Constraint Learning}
\ARTICLEAUTHORS{%
\AUTHOR{Yilin Zhu}
\AFF{Department of Mathematics, University of Iowa, \EMAIL{yilin-zhu@uiowa.edu}}
\AUTHOR{Samuel Burer}
\AFF{Department of Business Analytics, University of Iowa, \EMAIL{samuel-burer@uiowa.edu}}
}

}{}

\ifthenelse{\boolean{ijoc}}{}{\usepackage{amsmath}}
\usepackage{amssymb}
\ifthenelse{\boolean{ijoc}}{}{\usepackage{amsthm}}

\usepackage{setspace}
\usepackage{graphicx}
\usepackage[margin=2.54cm]{geometry}
\usepackage[pdftex,colorlinks,linkcolor=blue,citecolor=blue,urlcolor=blue]{hyperref}
\usepackage{mathtools}
\usepackage{siunitx}
\usepackage{multirow}

\renewcommand{\Box}{\text{\sc Box}}

\DeclareMathOperator{\st}{s. \hspace*{-0.03in} t.}

\DeclareMathOperator{\conv}{conv}

\DeclareMathOperator{\proj}{proj}

\DeclareMathOperator{\CH}{CH}

\ifthenelse{\boolean{ijoc}}{}
{
\newtheorem{proposition}{Proposition}%[section]
}

\title{An Extended Validity Domain for Constraint Learning}

\author{%
Yilin Zhu\thanks{Department of Mathematics, University
of Iowa, Iowa City, IA, 52242--1994, USA\@. Email: {\tt
yilin-zhu@uiowa.edu}.}%
\and
Samuel Burer\thanks{Department of Business Analytics, University
of Iowa, Iowa City, IA, 52242--1994, USA\@. Email: {\tt
samuel-burer@uiowa.edu}.}%
}

\date{May 27, 2025}

\begin{document}

\def\myabstract{
\noindent We consider embedding a predictive machine-learning model within a prescriptive optimization problem. In this setting, called {\em constraint learning\/}, we study the concept of a {\em validity domain\/}, i.e., a constraint added to the feasible set, which keeps the optimization close to the training data, thus helping to ensure that the computed optimal solution exhibits less prediction error. In particular, we propose a new validity domain which uses a standard convex-hull idea but in an extended space. We investigate its properties and compare it empirically with existing validity domains on a set of test problems for which the ground truth is known. Results show that our extended convex hull routinely outperforms existing validity domains, especially in terms of the function value error, that is, it exhibits closer agreement between the true function value and the predicted function value at the computed optimal solution. We also consider our approach within two stylized optimization models, which show that our method reduces feasibility error, as well as a real-world pricing case study.
}

\ifthenelse{\boolean{ijoc}}{\ABSTRACT{\myabstract}}{}

\maketitle

\ifthenelse{\boolean{ijoc}}{}%
{
\begin{abstract}
\myabstract
\end{abstract}
}

%\ifthenelse{\boolean{ijoc}}{}{\begin{onehalfspace}}

\section{Introduction} \label{sec:intro}

The fields of optimization and machine learning (ML) are closely
intertwined, a prominent example being the use of optimization as a
subroutine for training ML models. ML assists optimization in a number
of interesting and beneficial ways, too. \cite{osti_10337537} classify
such approaches into two groups: {\em ML-augmented optimization\/},
which uses ML to enhance existing optimization algorithms, e.g., when
ML models emulate expensive branching rules within a branch-and-bound
algorithm; and {\em end-to-end optimization\/}, which combines ML and
optimization techniques to obtain the optimal solutions of optimization
problems. \cite{De_Filippo_2018} and \cite{bengio2021machine} also
survey the use of ML techniques for modeling various components of
combinatorial-optimization algorithms, thus improving the algorithms'
accuracy and efficiency.

\cite{SADANA2025271} have recently proposed another taxonomy to
describe the interactions between ML and optimization. They focus on
{\em contextual optimization\/}, i.e., when an optimization problem
depends on uncertain parameters that are themselves correlated with
available side information, covariates, and features. In particular,
the authors identify three subcategories of contextual optimization:
(i) {\em decision-rule optimization\/} uses an ML model to predict
an optimal solution directly; (ii) {\em sequential learning and
optimization\/} first uses an ML model to predict uncertain parameters
in an optimization problem, which is subsequently solved to obtain an
optimal solution; and (iii) {\em integrated learning and optimization\/}
combines both training and optimization with the goal of improving the
quality of the final optimal solution, not specifically the quality of
the ML prediction.

In this paper, we consider a specific case of sequential
learning and optimization called {\em constraint learning
(CL)\/}; see \cite{maragno2025mixed} and the survey of
\cite{fajemisin2023optimization}. In this setting, a predictive model
is included as a function---which we denote by $\hat h$---within an
optimization problem. This function is learned from empirical data,
for example, because it lacks an explicit formula and hence cannot be
employed within an optimization model in a traditional manner. After
$\hat h$ is learned and inserted into a constraint or objective, the
optimization model is solved to obtain an optimal solution.

Because many popular ML models (e.g., linear regression, neural
networks, and decision trees) are mixed-integer-programming (MIP)
representable, the final optimization model can often be solved by
off-the-shelf software such as Gurobi. To simplify the implementation
and solution of CL models as MIPs, several software packages have
recently been developed. These include OptiCL \citep{maragno2025mixed},
JANOS \citep{bergman2019janos}, OMLT \citep{ceccon2022omlt}, PyEPO
\citep{tang2024pyepo}, and Gurobi versions 10 and later \citep{gurobi}.
While the precise details of how these packages embed an ML model into
an optimization problem are critical for the overall efficiency of the
solver's performance, we do not consider such details in this paper.

The CL paradigm has attracted significant attention recently.
\cite{tjeng2018evaluating} reformulated neural networks as MIPs to
evaluate their adversarial accuracy. \cite{maragno2025mixed} learned
a so-called {\em palatability constraint\/} for the optimization of
food baskets provided by the World Food Programme. The same authors
also optimized chemotherapy regimens for gastric-cancer patients. This
involved learning and embedding a clinically relevant total-toxicity
function within a constraint. In a hypothetical context related to
university admissions, \cite{bergman2019janos} maximized the expected
incoming class size using an ML model that predicts the probability
of individual students accepting an admission offer. Each student's
probability was a function of his or her high-school background as well
as a university decision variable for the amount of scholarship offered
to that student. The university also faced a fixed overall scholarship
budget. \cite{mistry2021mixed} used a gradient-boosted tree to model
the relationship between the proportions and properties of ingredients
within a concrete mixture, the goal being to optimize the strength of
such a mixture.

One challenge for CL is that the error inherent in the embedded model $\hat h$, as described above, may manifest as error in the final optimization model. Indeed, researchers have realized that the CL approach can sometimes lead to an unreasonable computed optimal solution $\hat x$. One particular downside occurs when $\hat x$, although feasible and optimal for the final optimization model based on $\hat h$, is nevertheless too far from the original data on which $\hat h$ has been trained. Because the accuracy of $\hat h$ can deteriorate far from the data (i.e., poor extrapolation), $\hat x$ can be optimal with respect to $\hat h$, while in reality being severely suboptimal---or even infeasible.

Researchers have proposed a remedy for this downside, called a {\em validity domain\/}. (Another common term in the literature is {\em trust region\/}, but since this term has already been used extensively in the nonlinear-programming literature, we prefer {\em validity domain\/} in this paper.) To guard against poor extrapolation and its downstream effect on the optimization, a validity domain further constrains the feasible region of the optimization to be closer to the data, i.e., to a subset where the predictions of $\hat h$ are   likely to be more reliable.

Many different types of validity domains have been proposed in the
literature. \cite{courrieu1994three} defined several validity domains
for the specific case when the learned function $\hat h$ is a neural
network. In particular, we will explore one of these in this paper: the
{\em convex hull\/} validity domain, i.e., when the variable $x$ of the
optimization is constrained to be within the convex hull of the training
data. \cite{schweidtmann2022obey} used persistent homology to study the
topological structure of the data and then constructed a validity domain
using the convex hull combined with one-class support vector machines.
\cite{maragno2025mixed} used an enlarged convex hull of the data set to
define a validity domain, thus allowing the final optimal solution to be
slightly outside the data. \cite{doi:10.1287/ijoc.2022.0312} compared
six different validity domain techniques, including one of their own
design based on an {\em isolation forest\/}. An isolation forest is
a type of one-class classification model that is, in particular,
MIP-representable. They showed that their isolation-forest validity
domain was generally the most accurate among the six, while still being
computationally efficient.

In this paper, we introduce a new validity domain, which is based on
applying the convex-hull idea in an extended space, which concatenates
the optimization variable $x$ with the output of the learned function
$\hat h$. This approach arises from an intuition that the convex-hull
idea is helpful not only for describing the original data set but
also for ``learning'' the constraints, objective, and optimal
solution of the optimization problem. We explain this intuition and
provide theoretical support in Section \ref{sec:extvdom}. In Sections
\ref{sec:numresults}--\ref{sec:case}, we test our approach and observe,
generally speaking, significant improvement in the error at the optimal
solution $\hat x$ of the final model. We also show that our approach
does not require significantly more time than the regular convex-hull
idea, which acts only in the $x$ space, while being faster than the
method of \cite{doi:10.1287/ijoc.2022.0312}.

Note that we focus on regression models instead of classification models; that is, the output of $\hat h$ is continuous, not discrete. Also, it is important to note that our approach is agnostic to the type of predictive model $\hat h$ used, i.e., our approach can be applied no matter the functional form of $\hat h$.

This paper is organized as follows. In Section \ref{sec:cl}, we recount required background on constraint learning, and then in Section \ref{sec:extvdom}, we introduce our new validity domain in the extended space. Sections \ref{sec:numresults}--\ref{sec:case} contain numerical results and examples illustrating our method. We conclude the paper in Section \ref{sec:conclusion}.

\section{Background on Constraint Learning} \label{sec:cl}

In this section, we provide the relevant background on constraint learning.

%The paper's main contribution is then introduced later in Section \ref{sec:extvdom}.

\subsection{Fundamentals} \label{sec:cl:fdmt}

%As discussed in the Introduction, {\em constraint learning (CL)\/} is a sequential-learning-and-optimization framework in which one or more predictive ML models are embedded within an optimization problem.

Formally, we study the following standard-form model introduced in \cite{fajemisin2023optimization}:
%\begin{subequations} \label{equ:stdform}
%\begin{align}
%    \hat v \quad := \quad \min_{x,y} \quad &f(x) \label{equ:stdform:f} \\
%    \st \quad \, &x \in X \label{equ:stdform:X} \\
%    &g(x) \le 0 \label{equ:stdform:g} \\ 
%    &\theta(y) \le 0 \label{equ:stdform:theta} \\
%    &y = \hat h(x). \label{equ:stdform:h}
%\end{align}
%\end{subequations}
%\begin{align}
%    \hat v \quad := \quad \min_{x,y} \quad &f(x) \nonumber \\
%    \st \quad \, &x \in X \nonumber \\
%    &g(x) \le 0 \label{equ:stdform} \tag{CL} \\ 
%    &\theta(y) \le 0 \nonumber \\
%    &y = \hat h(x). \nonumber
%\end{align}
\begin{equation} \label{equ:stdform}
    \hat v := \min_{(x,y) \in \widehat{F}} f(x)
    \quad \text{where} \quad 
    \widehat{F} := \{ (x,y) : x \in X, \ g(x) \le 0, \ \theta(y) \le 0, \ y = \hat h(x) \}.
\end{equation}
Here, $x \in \mathbb{R}^{n_1}$ is the vector of decision variables, $f : \mathbb{R}^{n_1} \to \mathbb{R}$ is a function capturing the objective function in $x$, $X$ is a simple ground set such as a box or a sphere, and $g : \mathbb{R}^{n_1} \to \mathbb{R}^{m_1}$ is a function capturing the constraints on $x$. In addition, $y \in \mathbb{R}^{n_2}$ is a vector of auxiliary variables, and $\theta : \mathbb{R}^{n_2} \to \mathbb{R}^{m_2}$ is a function expressing constraints on $y$. Finally, $\hat h : \mathbb{R}^{n_1} \to \mathbb{R}^{n_2}$ is a function, corresponding to a predictive ML model, which maps the decision variables $x$ to the auxiliary variables $y$. In words, $y$ are the predictions of $x$ via the learned function $\hat h$. In total, there are $n := n_1 + n_2$ variables, $m := m_1 + m_2$ inequality constraints, and $n_2$ equality constraints. %For convenience, we define
%\[
%\widehat{F} := \{ (x,y) : \text{(\ref{equ:stdform:X})--(\ref{equ:stdform:h})} \}
%\widehat{F} := \{ (x,y) : x \in X, \ g(x) \le 0, \ \theta(y) \le 0, \ y = \hat h(x) \}
%\]
%to be the feasible set of (\ref{equ:stdform}).
We assume that $\widehat{F} \ne \emptyset$ and that (\ref{equ:stdform}) attains its optimal value $\hat v$, and we use the notation $\widehat{\text{Opt}}$ to denote the optimal solution set.

%If desired, the variable $y$ can be eliminated from (\ref{equ:stdform}) by substituting $\hat h(x)$ for $y$ throughout. However,

The function $\hat h$ plays an important modeling role by acting as a surrogate, or approximation, for a true function $h : \mathbb{R}^{n_1} \to \mathbb{R}^{n_2}$, which is not known explicitly but can be learned via a suitable ML technique applied to empirical data. Accordingly, we consider the {\em true optimization model\/}
%\begin{equation} \label{equ:true} \ifthenelse{\boolean{ijoc}}{}{\tag{{\sc True}}}
%    v^* \quad := \quad \min_{x,y} \left\{ \ f(x) \ : (x,y) \in F \right\},
%\end{equation}
%where
%\[
%    F := \left\{ (x,y): x \in X, \ \ g(x) \le 0, \ \theta(y) \le 0, \ y = h(x) \ \right\}
%\]
% is the {\em true feasible set\/}. 
\begin{equation} \label{equ:true}
    v^* := \min_{(x,y) \in F} f(x)
    \quad \text{where} \quad 
    F := \{ (x,y) : x \in X, \ g(x) \le 0, \ \theta(y) \le 0, \ y = h(x) \},
\end{equation}
which differs from (\ref{equ:stdform}) only by the equation $y = h(x)$ in the constraints. The value $v^*$ denotes the {\em true optimal value\/}, and we use $\text{Opt}^*$ to denote the {\em true optimal solution set\/}, which we assume to be nonempty. Roughly speaking, if $\hat h$ is a good approximation of $h$, then one expects $\hat v$ and $\widehat{\text{Opt}}$ to be good approximations of $v^*$ and $\text{Opt}^*$.

We assume that the approximation $\hat h$ of $h$ is learned from a data set in $X$ of size $N$, 
\begin{equation} \label{equ:dataset}
D_N := \{ x^{(i)} \in X : i = 1, \ldots, N \},
\end{equation}
and from the corresponding $N$ observed function values
$h(D_N) := \left \{ h(x^{(i)}) : i = 1, \ldots, N \right\}$.
In practice, the observed function values contain noise, and ML techniques that learn $\hat h$ from $D_N$ and $h(D_N)$ should account for this noise, e.g., to avoid overfitting. In this paper, we will not explicitly model the noise. Rather, just as ML techniques learn $\hat h$ using the noisy evaluations of $h$, our method described in Section \ref{sec:extvdom} will make direct use of the noisy evaluations as well when applied in practice; see Sections \ref{sec:numresults}-\ref{sec:case}.

Other variations of (\ref{equ:stdform}) are possible. For example, the objective and constraint functions $f$ and $g$ could involve both $x$ and $y$, not just $x$. In addition, there could be auxiliary features, say $w$, incorporating contextual information about the optimization problem. Another possibility is additional variables, say $z$, that are (nonlinear) combinations of the main decision variables $x$. Together, $w$ and $z$ could then be used to build a better approximation $\hat h(x,w,z)$ of the predicted outputs. For the sake of simplicity, we will not incorporate $w$ and $z$ in this paper, but we refer the reader to \cite{fajemisin2023optimization} for additional references, which do take into account $w$ and $z$.

From time to time in this paper, we will also refer to an even simpler form of (\ref{equ:stdform}), which optimizes a predictive model over a simplified feasible set:
\begin{equation} \label{equ:justobj}
    \min_x \{ \ \hat f(x) \ : \ x \in X \ \}.
\end{equation}
This problem is an instance of (\ref{equ:stdform}) but with $y$ appearing in the objective: $\min_{x,y} \{ y : x \in X, y = \hat f(x) \}$ and $g$ and $\theta$ nonexistent.

\subsection{Errors} \label{sec:errors}

When comparing (\ref{equ:stdform}) with (\ref{equ:true}), hopefully $\hat v$ and $\widehat{\text{Opt}}$ are respectively close to $v^*$ and $\text{Opt}^*$. To measure this precisely, we formally define the errors
\begin{align*}
\text{optimal value error} &:= | \hat v - v^* |, \\
\text{optimal solution error} &:= \| \hat x - x^* \|,
\end{align*}
where $\hat x$ and $x^*$ are given optimal solutions for (\ref{equ:stdform}) and (\ref{equ:true}), respectively. Note that the optimal solution error depends on a {\em particular\/} pair of optimal solutions and hence could vary for different pairs if there are indeed multiple optimal solutions, but we choose this definition since it is both natural and practical.

Another type of error refers directly to the quality of the approximation $\hat h$ of $h$ for a particular feasible solution $(\hat x, \hat y)$ of (\ref{equ:stdform}). We define the
\[
\text{function value error at } (\hat x, \hat y) \ := \ \| \hat y - h(\hat x) \| = \| (\hat h - h)(\hat x) \|.
\]
That is, the function value error at $\hat x$ measures the Euclidean difference in $\mathbb{R}^{n_2}$ between the true value of $h$ at $\hat x$ and its predicted value under $\hat h$. Although the function value error is defined for any $(\hat x, \hat y)$, we will typically measure it at an optimal solution of (\ref{equ:stdform}).

A final type of error measures how close to true feasibility a given feasible $(\hat x, \hat y) \in \widehat{F}$ is:
\[
\text{feasibility error at } (\hat x, \hat y) \ := \| \max\{0, \theta(h(\hat x))\} \|,
\]
where the maximum between 0 and the vector $\theta(h(\hat x))$ is taken component-wise. Recall that $\widehat F$ and the true feasible set $F$ differ only in that $\widehat F$ uses the constraint $y = \hat h(x)$, while $F$ uses $y = h(x)$. In particular, given $(\hat x, \hat y) \in \widehat F$, the naturally corresponding element of $F$ is $(\hat x, h(\hat x))$. Because $\hat x$ already satisfies $\hat x \in X$ and $g(\hat x) \le 0$, true feasibility of $(\hat x, h(\hat x))$ then holds if and only if $\theta(h(\hat x)) \le 0$. The feasibility error is simply a measure of how much this constraint is violated.

We remark that both the function value and feasibility errors are exactly measurable only when $h$ can be evaluated without noise.

%\subsection{Implementation details}

%The practical application of CL requires that problem (\ref{equ:stdform}) be readily solvable to global optimality. This can depend, for example, on whether (\ref{equ:stdform}), and in particular $\hat h$, can be represented in a form that can be handled by one's preferred optimization solver. Fortunately, many standard predictive models give rise to $\hat h$, which can be represented using mixed-integer-programming (MIP) techniques. These include linear regression, decision trees, and neural networks as mentioned in the Introduction. To illustrate briefly, the ReLU activation function $y := \max\{0, x\}$ can be modeled using the MIP-system
%\begin{equation*} \label{equ:relu}
%    \begin{array}{ll@{}ll}
%    & y \geq x\\
%    & y \leq x+M(1-z) \\
%    & y \leq Mz \\
%    & y \geq 0\\
%    & z \in \{0,1\},
%    \end{array}
%\end{equation*}
%where $z$ is an extra auxiliary variable and $M > 0$ is a valid upper bound on the absolute value of any feasible $x$.

\subsection{Validity domains} \label{sec:cl:vdomain}

Solving (\ref{equ:stdform}) may lead to large errors compared to the true optimization problem (\ref{equ:true}). The {\em validity domain\/} concept has thus been introduced to mitigate these errors, particularly the function value error. Given a validity domain $V \subseteq \mathbb{R}^{n_1}$, we introduce the following optimization problem, which is closely related to (\ref{equ:stdform}):
%\begin{subequations} \label{equ:vdomform}
%\begin{align}
%    \hat v(V) \quad := \quad \min_{x,y} \quad &f(x) \label{equ:vdomform:f} \\
%    \st \quad \, &(x,y) \in \widehat{F} \label{equ:vdomform:Fhat} \\
%                 &x \in V
%\end{align}
%\end{subequations}
\begin{equation} \label{equ:vdomform}
    \hat v(V) := \min_{x,y} \left\{ f(x) :  (x,y) \in \widehat{F}, x \in V \right\}
\end{equation}
In words, (\ref{equ:vdomform}) is simply (\ref{equ:stdform}) with the added constraint $x \in V$. We denote the optimal value as $\hat v(V)$ to reflect the presence of $V$, and $\hat x(V)$ denotes an optimal solution of (\ref{equ:vdomform}).

We next describe several validity domains, which have been proposed in the literature.

\subsubsection{Simple bounds}

A simple validity domain is the smallest coordinate aligned hypercube, which surrounds the empirical data $D_N$ defined in (\ref{equ:dataset}):
\begin{equation} \label{equ:box}
    \Box := \left\{ x \in \mathbb{R}^{n_1} : \min_{i=1}^N x^{(i)}_j \le x_j \le \max_{i=1}^N x^{(i)}_j \ \ \forall \ j = 1,\ldots,n_1 \right\}.
\end{equation}

\subsubsection{The convex hull}

Recall that $D_N$ in (\ref{equ:dataset}) denotes the $N$ points in
$X$ over which $\hat h$ has been learned using the (possibly noisy)
function values $h(D_N)$. Let $\conv (D_N)$ be the smallest convex
subset containing $D_N$, i.e., its convex hull. \cite{courrieu1994three}
first used the convex hull as a validity domain in the case that $\hat
h$ modeled feed-forward neural networks, and \cite{maragno2025mixed}
also advocated the use of $\conv(D_N)$ as a validity domain. Throughout
the rest of this paper, we will use the notation
\begin{equation} \label{equ:CH}
\CH := \conv(D_N) 
\end{equation}
to denote this specific validity domain. Including the constraint $x
\in \CH$ in (\ref{equ:stdform}) amounts to adding extra variables
and linear constraints, which can be expensive when $N$ is large.
\cite{maragno2025mixed} discuss strategies to ease the computational
burden, e.g., column generation. They also proposed several variants of
CH, which can be implemented in modern optimization solvers.

%For example, in situations where the data set $D_N$ decomposes into clusters $D_N = \cup_{k = 1}^\ell D_{N_k}$, they defined a validity domain equal to the union of the convex hulls of each cluster, i.e., $\cup_{k=1}^\ell \conv( D_{N_k})$. This validity domain is MIP-representable using standard disjunctive-programming techniques that involve one binary variable for each cluster. They also introduced another validity domain, which slightly enlarges $\CH$ via a user-defined parameter $\varepsilon > 0$. In particular, they proposed the Minkowski sum $\CH + B(0, \varepsilon)$, where $B(0, \varepsilon)$ is the Euclidean ball of radius $\varepsilon$ centered at 0. Using second-order-cone constraints to represent $B(0,\varepsilon)$, this domain is straightforward to include in many modern optimization solvers.

\subsubsection{The enlarged convex hull} \label{sec:enlarged_ch}

\cite{balestriero2021learninghighdimensionamounts} argue that, for
higher-dimensional datasets, say, with $n \ge 100$, the phenomenon
of interpolation (i.e., when a new sample belongs to the convex
hull of the existing samples) almost surely never occurs. Hence,
generalization of a fitted model must occur via extrapolation, not
interpolation. As a result, in the context of CL, the convex hull method
just mentioned may be too conservative. To alleviate this conservatism,
\cite{maragno2025mixed} proposed an {\em enlarged convex hull\/}
formally defined as the Minkowski sum
\begin{equation} \label{equ:epsilon-CH}
\CH_{\varepsilon} := \CH + B_p(0,\varepsilon),
\end{equation}
where $B(0,\varepsilon)$ is the $p$-norm ball of a user-supplied radius
$\varepsilon > 0$ centered at 0. When $p=1$ or $p=\infty$, this can be
modeled using linear programming, or when $p=2$, using second-order
cone programming. The authors also demonstrated the effectiveness of
$\CH_\varepsilon$ in their setting.

\subsubsection{Support vector machines}

%A one-class support vector machine (OCSVM) is an unsupervised classification algorithm used for anomaly detection \cite{scholkopf1999support}, which works by defining the boundary between in-distribution and out-of-distribution data. Like many SVM approaches, it uses the kernel trick to implicitly map data into a higher-dimensional feature space so that determining the decision boundary is easier. In particular, a kernel function $K(\cdot,\cdot)$ is used to calculate the inner product between input vectors in the feature space, and several different kernel functions such as linear, polynomial, sigmoidal, and Gaussian (also known as RBF) have been employed \cite{vapnik2013nature}. The boundary is then defined in the variable $x$ by the equation $\sum_{i} \alpha_iK(x,x_i)+b=0$, where $i$ indexes the support vectors, $x_i$ are the support vectors, and $\alpha_i, b$ are the learned parameters. In particular, a new point $x$ with $\sum_{i} \alpha_iK(x,x_i)+b \ge 0$ is classified as in-distribution.

\cite{schweidtmann2022obey} proposed the use of a one-class support vector machine to define a validity domain. In particular, the authors expressed $x$ as being in-distribution using a validity domain defined by the constraint
$
    \sum_{i}\alpha_i K \left(x_i, x \right) + b \geq 0,
$
where $i$ indexes the support vectors, $x_i$ are the support vectors, $K(x_i,x)$ is the RBF kernel function $\exp \left( -\gamma \cdot \|x-x_i\|^2 \right)$ with a hyperparameter $\gamma$, and $\alpha_i, b$ are the learned parameters. Although this constraint is neither convex nor MIP-representable, the authors demonstrated numerical success using a custom global optimization algorithm.

%A similar, albeit MIP-representable technique was introduced by \cite{shang2017data}. They proposed a piecewise linear kernel called the generalized intersection kernel and showed convexity of the corresponding MIP formulation.

\subsubsection{Isolation forests}

Similar to the preceding subsection, another model for outlier detection is the {\em isolation forest\/} \citep{liu2008isolation}, which is an ensemble-tree method to isolate outliers based on their quantitative characteristics, e.g., their atypical feature values. Each tree in the forest is trained to isolate instances by their features, and by design, outliers are closer to the root of the tree because they are more susceptible to isolation. The isolation forest then defines an instance to be an outlier if it has a short average path length to the root node, where the average is taken over all isolation trees in the forest.

Using a given isolation forest trained on the observed input data set,
\cite{doi:10.1287/ijoc.2022.0312} defined a validity domain to be the
set of points $x$, where the path length of $x$ in every tree of the
forest is larger than a pre-determined threshold $d$. Note that this is
equivalent to disallowing any $x$ that is an outlier in some tree, thus
slightly deviating from the use of the average measure described in the
preceding paragraph. They also developed a MIP representation of the
isolation forest enabling it to be embedded into an optimization model
and showed the effectiveness of this validity domain compared to others
in the literature.

\section{An Extended Validity Domain} \label{sec:extvdom}

We introduce a new validity domain, which is based on a simple observation---though perhaps not well known---that (\ref{equ:true}) is equivalent to a convex optimization problem in an extended space. Note that, as discussed in Section \ref{sec:cl}, we assume throughout this section that noiseless observations of the true function $h(x)$ are available, even as the precise functional form of $h(x)$ is unknown.

\subsection{Intuition}

Recall that $F$ denotes the true feasible region of (\ref{equ:true}) in the space $(x,y) \in \mathbb{R}^n$. Now define $$F^+ := \{ (x, y, f(x)) : (x,y) \in F \}$$ to be the {\em extended feasible set\/}, which is also sometimes called the {\em graph of $f$ over $F$\/}. In words, $F^+$ is the set of all triples $(x,y,\phi)$ in $\mathbb{R}^{n + 1}$, where $(x,y)$ is feasible and the scalar $\phi = f(x)$ encodes the objective value. The following proposition states that (\ref{equ:true}) is equivalent to the minimization of $\phi$ over $(x,y,\phi) \in \conv(F^+)$, i.e., a convex minimization:

\begin{proposition} \label{pro:equivconvprob}
$v^* = \min\{ \phi : (x,y,\phi) \in \conv(F^+) \}$.
\end{proposition}

\def\myprooftext1{%
Note that $\phi$ is a scalar variable, which represents the function value $f(x)$ in the extended space. Let $\nu$ denote the optimal value of $\min\{ \phi : (x,y,\phi) \in \conv(F^+) \}$, and let $(x^*,y^*)$ denote an optimal solution of (\ref{equ:true}) such that $f(x^*) = v^*$. Then $(x^*,y^*,v^*)$ is a member of $\conv(F^+)$, and hence $v^* \ge \nu$. To show $v^* \le \nu$, consider an optimal $(x,y,\nu) \in \conv(F^+)$. There exists a finite index set $K \ni k$, feasible points $(x^k,y^k)$ for (\ref{equ:true}), and multipliers $\lambda_k \ge 0$ such that $\sum_{k \in K} \lambda_k = 1$ and
\[
(x,y,\nu) = \sum_{k \in K} \lambda_k (x^k, y^k, f(x^k)) = \left( \sum_{k \in K} \lambda_k x^k,
\sum_{k \in K} \lambda_k y^k, \sum_{k \in K} \lambda_k f(x^k) \right).
\]
In particular, $\nu$ is greater than or equal to the minimum $f(x^k)$, which is itself greater than or equal to $v^*$.
}

\ifthenelse{\boolean{ijoc}}{
\proof{Proof.}
\myprooftext1
$\Halmos$
}{
\begin{proof}
\myprooftext1
\end{proof}
}

Proposition \ref{pro:equivconvprob} establishes that solving (\ref{equ:true}) is equivalent to the minimization of a linear function over the convex hull of $F^+$. Hence, in a certain sense, understanding $\conv(F^+)$ provides a key to solving (\ref{equ:true}). Indeed, in what follows, we take the point of view that $\conv(F^+)$ can help us ``learn'' the optimization problem (\ref{equ:true}) from a new perspective---one which is complementary to the use of $\hat h$ described in Section \ref{sec:cl}.

Proposition \ref{pro:inthelimit} below provides further intuition that $\conv(F^+)$ learns and describes the true optimization model (\ref{equ:true})---but this time from a data-driven perspective. Let $D_N$ be the data set given in (\ref{equ:dataset}) over which $\hat h$ has been trained. We define two related data sets:
\begin{align}
    D_N^+ &:= \left\{ \left( x^{(i)}, y^{(i)}, \phi^{(i)} \right) \ : \
        y^{(i)} = h(x^{(i)}), \  \phi^{(i)} = f(x^{(i)}), \ i = 1, \ldots, N \right\} \label{equ:dataset+} \\
    F_N^+ &:= \left\{ (x,y,\phi) \in D_N^+ : (x,y) \in F \right\} \label{equ:feasdataset+}
\end{align}
Here, $D_N^+ \subseteq \mathbb{R}^{n + 1}$ is an extension of $D_N$, which appends the true values of $h$ and $f$ to each data point $x^{(i)}$. Further, $F_N^+$ is the restriction of $D_N^+$ to just those data points that are feasible for (\ref{equ:true}). In particular, $D_N^+$ may contain both feasible and infeasible $(x^{(i)}, y^{(i)})$. For these definitions, recall our assumption that noiseless evaluations of $h$ are available, which in particular allows one to check which data points in $D_N^+$ are feasible.

Using the extended data sets $D_N^+$ and $F_N^+$ and viewing the sample size $N$ as a parameter, we now show that, if one takes larger and larger samples such that $F_N^+$ becomes dense in the extended feasible set $F^+$, then solving the surrogate optimization over $\conv(F_N^+)$ eventually solves (\ref{equ:true}).

\begin{proposition} \label{pro:inthelimit}
Let $\{ D_N \}_{N=1}^\infty \subseteq X$ be a sequence of data sets with $$\conv(F^+) \subseteq \lim_{N \to \infty} \conv(F_N^+),$$ and define $v_N^* := \min \{ \phi : (x,y,\phi) \in \conv(F_N^+) \}$. Then $v^* = \lim_{N \to \infty} v_N^*$.
\end{proposition}

\def\myprooftext2{%
Because $\conv(F_N^+) \subseteq \conv(F^+)$ for all $N$, it holds that $v_N^* \ge v^*$ for all $N$. Now consider $(x^*, y^*, \phi^*) \in F^+$, where $(x^*, y^*)$ is an optimal solution of (\ref{equ:true}) and hence $\phi^* = v^*$. By assumption, there exists a sequence $\{ (x_N, y_N, \phi_N) \in \conv(F_N^+) \}$ converging to $(x^*, y^*, v^*)$. In particular, $\{ \phi_N \} \to v^*$. Note also that $\phi_N \ge v^*_N$ since $(x_N, y_N, \phi_N)$ is a member of $\conv(F_N^+)$. In total, $\phi_N \ge v^*_N \ge v^*$ for all $N$ with $\{ \phi_N \} \to v^*$.  Hence, $\{ v_N^* \} \to v^*$.
}

\ifthenelse{\boolean{ijoc}}{
\proof{Proof.}
\myprooftext2
$\Halmos$
}{
\begin{proof}
\myprooftext2
\end{proof}
}

\subsection{Our new validity domain and its variants} \label{sec:new_vd_and_variants}

Based on the prior subsection, we propose the following validity domain, where $F_N^+$ is the extended data set of feasible points defined in (\ref{equ:feasdataset+}) based on $D_N$ and $D_N^+$ from (\ref{equ:dataset}) and (\ref{equ:dataset+}), respectively:
\begin{equation} \label{equ:CH+}
\CH^+ := \conv(F_N^+).
\end{equation}

\noindent While $\CH \subseteq \mathbb{R}^{n_1}$ from (\ref{equ:CH}) is defined in the space of $x$, this new validity domain $\CH^+ \subseteq \mathbb{R}^{n+1}$ is defined in the space of $(x,y,\phi)$. However, it is easy to see that the projection of $\CH^+$ onto the variable $x$ is contained in $\CH$, i.e., $\proj_x(\CH^+) \subseteq \CH$.

In order to use $\CH^+$ as a validity domain, we simply add the constraint $(x,y,f(x)) \in \CH^+$ to (\ref{equ:stdform}):
%\begin{align*}
%    \hat v(\CH^+) \quad := \quad \min_{x,y} \quad &f(x) \\
%    \st \quad \, &(x,y) \in \widehat{F} \\
%                 &(x,y,f(x)) \in \CH^+.
%\end{align*}
\[
    \hat v(\CH^+) := \min_{x,y} \left\{ f(x) :
    (x,y) \in \widehat{F}, (x,y,f(x)) \in \CH^+ \right\}.
\]
Note that $\hat h$ appears in both constraints through the equation $y = \hat h(x)$. Based on the fact that $\proj_x(\CH^+) \subseteq \CH$, we have the following immediate relationship between $\hat v(\CH^+)$ and $\hat v(\CH)$.

\begin{proposition}
$\hat v(\CH^+) \ge \hat v(\CH)$.
\end{proposition}

In fact, we believe the property $\proj_x(\CH^+) \subseteq \CH$ is a defining feature of our approach. The validity domain $\CH$ acts as a natural geometric restriction on $x$ to keep the optimization close to the training data---with a goal to ameliorate the function value error. By further constraining $\CH$ in the space $x$, our intuition is that the new validity domain $\CH^+$ further aids the optimization as argued in Propositions \ref{pro:equivconvprob}--\ref{pro:inthelimit}. Practically speaking, in Sections \ref{sec:numresults}--\ref{sec:case}, we will show via example that $\CH^+$ does indeed further reduce the function value error relative to $\CH$.

Other variations of $\CH^+$, which maintain the property that the
projection onto $x$ is contained in $\CH$, are possible,
and such variations can be useful in practice (see Section
\ref{sec:case}). For example, consider a case in which
the historical data $F_N^+$ employed in the definition of $\CH^+$
in (\ref{equ:CH+}) is incomplete in the sense that the present-day
optimization problem includes additional decision variables for which
past data has not been collected. Said differently, one may have historical data on only a subset of the decision variables. Another example occurs when, even if
the past data is complete, a new constraint is imposed in the present,
rendering much of the historical data infeasible and hence theoretically
unfit for the construction of $\CH^+$. In such a case, it may prove
beneficial to make use of all past data---both feasible and infeasible
points---so as to not discard too much data.

More formally, let us consider the case that we do indeed have historical data $D_N$ on all variables $x$, and define $F_N := D_N \cap F$ to be those data points that are feasible for (\ref{equ:true}). Given arbitrary subsets $J \subseteq \{1,\ldots,n_2\}$
and $D \subseteq D_N$ such that $F_N \subseteq D$, consider the following variant of $\CH^+$:
\[
    \CH' := \conv\left( \{ (x,y_J, f(x)) : x \in D, y_J = h_J(x) \} \right)\subseteq \mathbb{R}^{n_1 + |J| + 1}.
\]
Note that $\CH'$ equals $\CH^+$ when $J = \{1,\ldots,n_2\}$ and $D = F_N$; see the definitions (\ref{equ:dataset+}) and
(\ref{equ:feasdataset+}). It is then clear that $\proj_x(\CH^+)
\subseteq \proj_x(\CH') \subseteq \CH$. In this sense, $\CH'$ is
sandwiched between $\CH^+$ and $\CH$. Another variant of $\CH^+$ could
drop the function value $f(x)$ in the definition of $\CH'$. In Sections \ref{sec:stylized}--\ref{sec:case}, we will show practical cases in which we choose a validity domain $\CH'$ between
$\CH$ and $\CH^+$. Even in those cases, we will call these the ``$\CH^+$
approach'' for simplicity.

Following Section \ref{sec:enlarged_ch}, we also consider the enlarged version of $\CH^+$, which is based on a user-supplied $\varepsilon > 0$:
\begin{equation} \label{equ:epsilon-CH+}
\CH^+_\varepsilon := \CH^+ + B_p(0, \varepsilon),
\end{equation}
where $B_p(0, \varepsilon)$ is the $p$-norm ball of radius $\varepsilon$ in $\mathbb{R}^{n+1}$. In our experiments in Section \ref{sec:numresults}, we will test $p=2$, which can be modeled using second-order programming.

\subsection{Illustration}

For illustration, we refer to the simplified form (\ref{equ:justobj}) described in Section \ref{sec:cl:fdmt}, which optimizes a true but unknown objective function $f(x)$ over a simple domain $x \in X$ by substituting a learned approximation $\hat f$ of $f$. Our problem in this subsection is thus $\min \{ \hat f(x) : x \in X \}$.

\begin{figure}[htbp]
\centering
\includegraphics[width=0.95\linewidth]{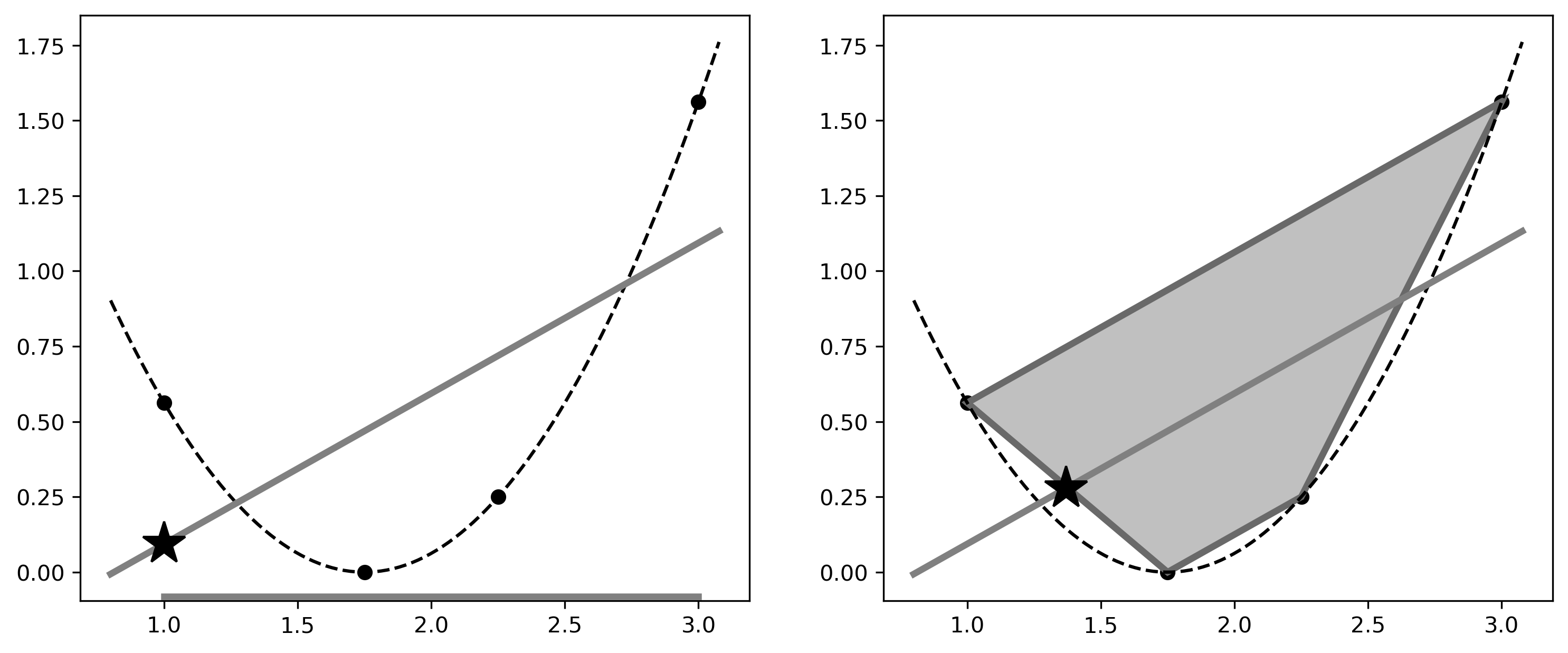}
\caption{Illustration of $\CH$ on the left and $\CH^+$ on the right. $\CH^+$ exhibits better function-value and optimal-solution errors.}
\label{fig:parabola_good}
\end{figure}

Figure \ref{fig:parabola_good} depicts a one-dimensional example ($n_1 = 1$) in which the true function to optimize is $f(x) = (x - 1.75)^2$ over all $x \in X := [0,4]$. The true optimal value is $v^* = 0$, and the true unique optimal solution is $x^* = 1.75$. We depict $f(x)$ as a dotted curve in both panels of Figure \ref{fig:parabola_good}, keeping in mind that this curve would be unknown in practice. We take $N = 4$ and $D_N = \{ 1.00, 1.75, 2.25, 3.00 \}$, and then $\hat f(x)$ is taken to be the least-squares regression line based on the true function $f$ evaluated at $D_N$ without noise; this is depicted as the solid, sloped line in both panels.

In the left panel, we also depict the validity domain $\CH$ as the horizontal line segment overlaid on the $x$-axis from $1$ to $3$, which is the convex hull of $D_N$. Then $\min \{ \hat f(x) : x \in X \cap \CH \}$ is optimized at $\hat x(\CH) = 1.00$. On the graph of $\hat f(x)$, this is depicted as the star in the left panel. In contrast, the right panel depicts the set $\CH^+$, which is the filled quadrilateral in $(x, y)$, spanning the four sampled points in $D_N^+$, which lie on the graph of $f(x)$. Then $\min \{ \hat f(x) : x \in X, \ (x, \hat f(x)) \in \CH^+ \}$ optimizes the function value over the intersection of the regression line and the filled quadrilateral. In this case, the optimal solution occurs at $\hat x(\CH^+) \approx 1.37$ and is depicted by the star in the right panel. The plots show that the approach based on $\CH^+$ exhibits better function-value and optimal-solution errors but worse optimal-value error.

Figure \ref{fig:parabola_good} also demonstrates that the $\CH^+$ method does not simply return the sampled point with minimum $f$ value. Indeed, the interaction of the sampled values $(x^{(i)}, f(x^{(i)}))$ with the function approximation $\hat f$ is critical to the behavior of $\CH^+$.

\section{Numerical Results} \label{sec:numresults}

To evaluate the extended validity domain $\CH^+$ introduced in
Section \ref{sec:extvdom}, we adopt a procedure introduced by
\cite{doi:10.1287/ijoc.2022.0312}. To this end, in the following
subsections, we define the concept of a basic experiment, then describe
our method for generating multiple experiments, and finally detail the
optimization results. All numerical experiments were coded using Python
3.10.8 and Gurobi 11.0 and conducted on a single Xeon E5-2680v4 core
running at 2.4 GHz with 16 GB of memory under the CentOS
Linux operating system. The code and results are publicly shared at
\url{https://github.com/yillzhu/extvdom}.

\subsection{Definition of an experiment}

In our testing, we define an {\em experiment\/} to be the full specification of six design options:

\begin{itemize}

    \item {\em Ground truth\/}, i.e., a true optimization problem of the form
    $
    v^* = \min_x \left\{ f(x) : x \in X \right\}
    $
    for which the function $f(x)$ and the true optimal value $v^*$ are known. A true optimal solution $x^* \in \text{Opt}^*$ is known as well. This corresponds to the simplified optimization (\ref{equ:justobj}) discussed in Section \ref{sec:cl:fdmt} in which the learned function lies in the objective, not the constraints. We will interchangeably identify a ground truth with its function $f$.

    \item {\em Sampling rule\/}, i.e., a rule ${\cal R}$ describing how to sample points in $X$.

    \item {\em Sample size\/}, i.e., the number of points $N$ to sample in $X$.

    \item {\em Noise level\/}, i.e., a scale factor $\sigma$ corresponding to the amount of random noise added to the function evaluations of $f$ during sampling.

    \item {\em Seed\/}, i.e., the seed $s$ used to initiate the random number generator before sampling.

    \item {\em ML technique\/}, i.e., the machine-learning technique ${\cal M}$ used to learn $\hat f$ from noisy evaluations of $f$.
    
\end{itemize}

Once the ground truth $f$, sampling rule ${\cal R}$, sample size $N$, noise factor $\sigma$, seed $s$, and ML model ${\cal M}$ are specified, an experiment proceeds as follows. Setting the seed $s$, we randomly sample $D_N \subseteq X$ following the sampling rule ${\cal R}$, and then we evaluate $f$ at all points in $D_N$ adding random noise scaled by the noise factor $\sigma$. Then we use the ML technique ${\cal M}$ to learn the approximation $\hat f$ using the empirical data $D_N$ and the noisy function values $f(D_N)$. The experiment continues by solving (\ref{equ:justobj}) based on $\hat f$ with an additional validity-domain constraint for several different domains.

In short, an experiment specifies the six design options, builds the optimization model based on $\hat f$, and  solves the model multiple times, each time testing a different validity domain. Because the ground truth is known, the errors defined in Section \ref{sec:errors} are easily computed for each validity domain. This allows us to determine, for a given experiment, which validity domain yields smaller errors. Finally, by aggregating these errors over multiple experiments, we can identify trends in the performance of different validity domains.

\subsection{Generating multiple experiments} \label{sec:numresults:setup}

Following \cite{doi:10.1287/ijoc.2022.0312}, we test seven different ground truths corresponding to seven challenging nonlinear benchmark functions from the literature (see, for example, \cite{testfunctions}), namely
$
f \in \left\{ \text{\em Beale, Griewank, Peaks, Powell, Qing, Quintic, Rastrigin} \right\}.
$
The input dimension $n_1$ of these functions varies from 2 to 10, and in each case, $X$ is an $n_1$-dimensional box. Further, we consider two sampling rules
$
{\cal R} \in \{ \text{\em Uniform, Normal} \},
$
where {\em Uniform\/} indicates a uniform sample over the box domain $X$ and {\em Normal\/} indicates a jointly independent normal sample around the global minimum $x^* \in X$ of $f$ with covariance matrix $\rho I$, where $I$ is identity matrix of size $n_1$ and $\rho$ is problem-specific. In particular, we choose $\rho$ to be $1/6$ of the distance of $x^*$ to the boundary of $X$. Due to the nature of the normal distribution in each dimension, this ensures samples $D_N$ that tend to be close to $x^*$ and highly unlikely to be outside $X$.

The sample size is taken as $N \in \bigl\{ 500,  1000,  1500 \bigr\}$, and the noise is a zero-mean univariate normal distribution with standard deviation equal to $\sigma$ times the standard deviation of the noisy $f(D_N)$ function values, where the scale factor $\sigma \in \bigl\{0.0, 0.1, 0.2 \bigr\}$. In particular, $\sigma = 0.0$ corresponds to the no-noise case. Furthermore, we take 100 seeds, specifically $s \in \bigl\{ 2023, 2024, \dots, 2122 \bigr\}$.

In addition, we test three different machine learning models ${\cal M}$ using the scikit-learn package of Python:
$
{\cal M} \in \left\{ \text{\em RandomForestRegressor, GradientBoostingRegressor, MLPRegressor\/} \right\}.
$
In particular, we implement {\em RandomForestRegressor\/} with 100 trees and maximum depth of 5; {\em GradientBoostingRegressor\/} with 100 boosting stages and a maximum depth of 5 for the individual regression estimators; and {\em MLPRegressor\/} with 2 hidden layers and 30 neurons in each layer. (For each experiment independently, we also tried using grid search on the parameters to find the best model for that experiment. Ultimately, we found that the overall conclusions about the various validity domains in Section \ref{sec:numresults:results} were very similar. So we have fixed the parameters of ${\cal M}$ to simplify our experiments and reduce testing time.) Before constructing the ML models, the input features are standardized in $[0,1]$ using min-max scaling, and the noisy function values are normalized to mean 0 and standard deviation 1.

We generate experiments by looping over all combinations of $(f,
{\cal R}, N, \sigma, s, {\cal M})$ for a total of $7 \times 2 \times
3 \times 3 \times 100 \times 3 = 37,800$ experiments. To provide a
snapshot of the quality of the ML models, as well as the time required
to compute them, Table \ref{tab:r2} shows the average $R^2$ scores
broken down by ground truth and ML technique, and similarly Table
\ref{tab:ml_train_times} shows the median training times (in seconds).
Specifically, for the calculation of each $R^2$ value, we train and test
the model using an 80-20 split of the data, and Table
\ref{tab:r2} shows the performance over the test set. However, the
models used for optimization as described in the next subsection are
trained on the complete data set.

% latex table generated in R 4.3.3 by xtable 1.8-4 package
% Fri May 31 08:20:53 2024
\begin{table}[htbp]
\centering%\small
\begin{tabular}{l|ccccccc}
 & Beale & Griewank & Peaks & Powell & Qing & Quintic & Rastrigin \\ 
  \hline
{\em RandomForestRegressor} & 0.88 & 0.34 & 0.87 & 0.55 & 0.66 & 0.82 & 0.17 \\ 
  {\em GradientBoostingRegressor} & 0.90 & 0.39 & 0.97 & 0.83 & 0.89 & 0.92 & 0.42 \\ 
  {\em MLPRegressor} & 0.93 & 0.62 & 0.90 & 0.91 & 0.89 & 0.95 & 0.04 \\ 
\end{tabular}
\caption{Mean test $R^2$ scores over all 37,800 experiments, grouped by function and ML technique.}
\label{tab:r2}
\end{table}

% latex table generated in R 4.3.3 by xtable 1.8-4 package
% Mon Jun  3 14:55:48 2024
\begin{table}[htbp]
\centering%\small
\begin{tabular}{l|ccccccc}
 & Beale & Griewank & Peaks & Powell & Qing & Quintic & Rastrigin \\ 
  \hline
{\em RandomForestRegressor} & 0.35 & 0.46 & 0.34 & 0.47 & 0.65 & 0.50 & 0.75 \\ 
  {\em GradientBoostingRegressor} & 0.24 & 0.41 & 0.24 & 0.41 & 0.72 & 0.49 & 0.89 \\ 
  {\em MLPRegressor} & 2.42 & 1.75 & 1.64 & 1.76 & 1.55 & 1.23 & 2.19 \\ 
\end{tabular}
\caption{Median training times (in seconds) over all 37,800 experiments, grouped by function and ML technique.}
\label{tab:ml_train_times}
\end{table}

\subsection{Optimization results for four validity domains} \label{sec:numresults:results}

For each experiment, we solve the corresponding optimization model four times by varying the validity domain 
$
V \in \left\{ \text{{\sc Box}, CH, {\sc IsoFor}, CH$^+$} \right\},
$
where {\sc Box}, CH, and CH$^+$ are defined in Sections \ref{sec:cl}--\ref{sec:extvdom}, and {\sc IsoFor} refers to the isolation-forest validity domain of \cite{doi:10.1287/ijoc.2022.0312}, which has also been described in Section \ref{sec:cl}. For {\sc IsoFor}, we follow \cite{doi:10.1287/ijoc.2022.0312} by setting hyperparameters to their default values and by taking the maximum depth of a tree to be 5 for the {\em Beale} and {\em Peaks} functions and 6 otherwise. In the online supplemental material, we show the median optimization setup and solve times for all combinations of ML technique and validity domains. 

We now examine the median errors associated with
each validity domain. Table \ref{tab:all} presents results for all
experiments grouped by function, type of error, and sampling rule. Each
group then has four errors corresponding to the four validity domains.
Furthermore, within each group of four errors, we scale so that the
median error corresponding to {\sc Box} equals 1.00, and
we show the smallest eror in a bold font, thus facilitating comparison
of the different validity domains. Note that we choose
to report medians rather than means to avoid any out-sized influence
by the tails of the associated distributions. However, in the online
supplemental material, we report the means of the same experiments for
comparison; we find that the qualitative conclusions do not change.

% latex table generated in R 4.3.3 by xtable 1.8-4 package
% Wed Jun 12 18:49:07 2024
\begin{table}[htbp]
\centering
\begin{tabular}{ll||cc|cc|cc}
Function &
\begin{tabular}[c]{@{}c@{}} Validity \\ Domain \end{tabular} &
\multicolumn{2}{c|}{%
\begin{tabular}[c]{@{}c@{}} Median Function \\ Value Error \end{tabular}%
} &
\multicolumn{2}{c|}{%
\begin{tabular}[c]{@{}c@{}} Median Optimal \\ Value Error \end{tabular}%
} &
\multicolumn{2}{c}{%
\begin{tabular}[c]{@{}c@{}} Median Optimal \\ Solution Error \end{tabular}%
} \\ 
& & {\em Uniform\/} & {\em Normal\/} & {\em Uniform\/} & {\em Normal\/} & {\em Uniform\/} & {\em Normal\/} \\
  \hline\hline
Beale & {\sc Box} & 1.00 & 1.00 & 1.00 & 1.00 & 1.00 & 1.00 \\ 
   & $\CH$ & 0.97 & 0.87 & 1.00 & 0.97 & 1.01 & 1.18 \\ 
   & {\sc IsoFor} & 0.63 & 0.76 & 0.73 & 0.91 & \textbf{0.81} & \textbf{0.52} \\ 
   & $\CH^+$ & \textbf{0.09} & \textbf{0.35} & \textbf{0.16} & \textbf{0.72} & 0.86 & 0.79 \\ \hline
  Griewank & {\sc Box} & 1.00 & 1.00 & \textbf{1.00} & 1.00 & 1.00 & 1.00 \\ 
   & $\CH$ & 0.73 & 0.95 & 1.09 & 1.00 & 0.92 & 0.97 \\ 
   & {\sc IsoFor} & 0.86 & 1.00 & 1.56 & 1.00 & \textbf{0.21} & \textbf{0.90} \\ 
   & $\CH^+$ & \textbf{0.49} & \textbf{0.53} & 1.18 & \textbf{0.90} & 0.89 & 1.08 \\ \hline
  Peaks & {\sc Box} & \textbf{1.00} & 1.00 & \textbf{1.00} & 1.00 & 1.00 & 1.00 \\ 
   & $\CH$ & \textbf{1.00} & 1.00 & \textbf{1.00} & 1.00 & 1.01 & 1.01 \\ 
   & {\sc IsoFor} & 1.47 & 1.25 & 1.16 & 1.01 & 1.04 & 0.88 \\ 
   & $\CH^+$ & 1.02 & \textbf{0.68} & \textbf{1.00} & \textbf{0.94} & \textbf{0.93} & \textbf{0.92} \\  \hline
  Powell & {\sc Box} & 1.00 & 1.00 & 1.00 & 1.00 & 1.00 & 1.00 \\ 
   & $\CH$ & 0.99 & 1.03 & 0.98 & 0.91 & 0.95 & 0.95 \\ 
   & {\sc IsoFor} & 0.90 & 1.06 & 0.89 & 0.79 & \textbf{0.66} & \textbf{0.59} \\ 
   & $\CH^+$ & \textbf{0.09} & \textbf{0.17} & \textbf{0.15} & \textbf{0.23} & 0.78 & 0.63 \\  \hline
  Qing & {\sc Box} & 1.00 & 1.00 & 1.00 & 1.00 & 1.00 & 1.00 \\ 
   & $\CH$ & 0.61 & 0.54 & 1.00 & 1.00 & 0.83 & 0.83 \\ 
   & {\sc IsoFor} & 1.20 & 0.73 & 1.00 & 1.00 & 0.75 & 0.72 \\ 
   & $\CH^+$ & \textbf{0.41} & \textbf{0.45} & \textbf{0.79} & \textbf{0.95} & \textbf{0.70} & \textbf{0.54} \\  \hline
  Quintic & {\sc Box} & 1.00 & 1.00 & 1.00 & 1.00 & 1.00 & 1.00 \\ 
   & $\CH$ & 0.63 & 0.42 & \textbf{0.90} & 0.36 & 0.97 & 0.97 \\ 
   & {\sc IsoFor} & 0.25 & 0.16 & 0.94 & 0.26 & \textbf{0.77} & 0.75 \\ 
   & $\CH^+$ & \textbf{0.13} & \textbf{0.06} & 1.00 & \textbf{0.25} & 0.92 & \textbf{0.74} \\  \hline
  Rastrigin & {\sc Box} & 1.00 & 1.00 & \textbf{1.00} & \textbf{1.00} & 1.00 & 1.00 \\ 
   & $\CH$ & 0.87 & 0.66 & 1.02 & 1.49 & 0.88 & 0.65 \\ 
   & {\sc IsoFor} & 0.92 & 0.84 & 1.04 & 1.39 & 0.82 & \textbf{0.63} \\ 
   & $\CH^+$ & \textbf{0.68} & \textbf{0.49 }& 1.13 & 1.54 & \textbf{0.70} & \textbf{0.63}
\end{tabular}
\caption{Errors for all 37,800 experiments, grouped by function, type of error, and sampling rule. In each group of four errors corresponding to four validity domains, the errors are scaled so that {\sc Box} has value 1.00 and the smallest error is shown in a bold font} in order to facilitate comparison among the different validity domains.
\label{tab:all}
\end{table}

For both sampling rules, we see small function value errors for $\CH^+$. In particular, for six of the seven functions, $\CH^+$ achieves the best median function value error for {\em Uniform\/}. For the seventh function ({\em Peaks\/}), $\CH^+$ achieves nearly the best---1.02 versus 1.00 for both {\sc Box} and $\CH$. For the {\em Normal\/} sampling rule, $\CH^+$ also performs well, where it achieves the best median function value error for all seven functions.

For the optimal value and optimal solution errors, $\CH^+$ performs
competitively. With respect to {\em Uniform\/}, for five of the seven
functions, $\CH^+$ has either the best median optimal value error or
the best median optimal solution error. For {\em Normal\/}, we see this
performance for all seven functions.

\begin{figure}[htbp]
\centering
\includegraphics[width=1.0\linewidth]{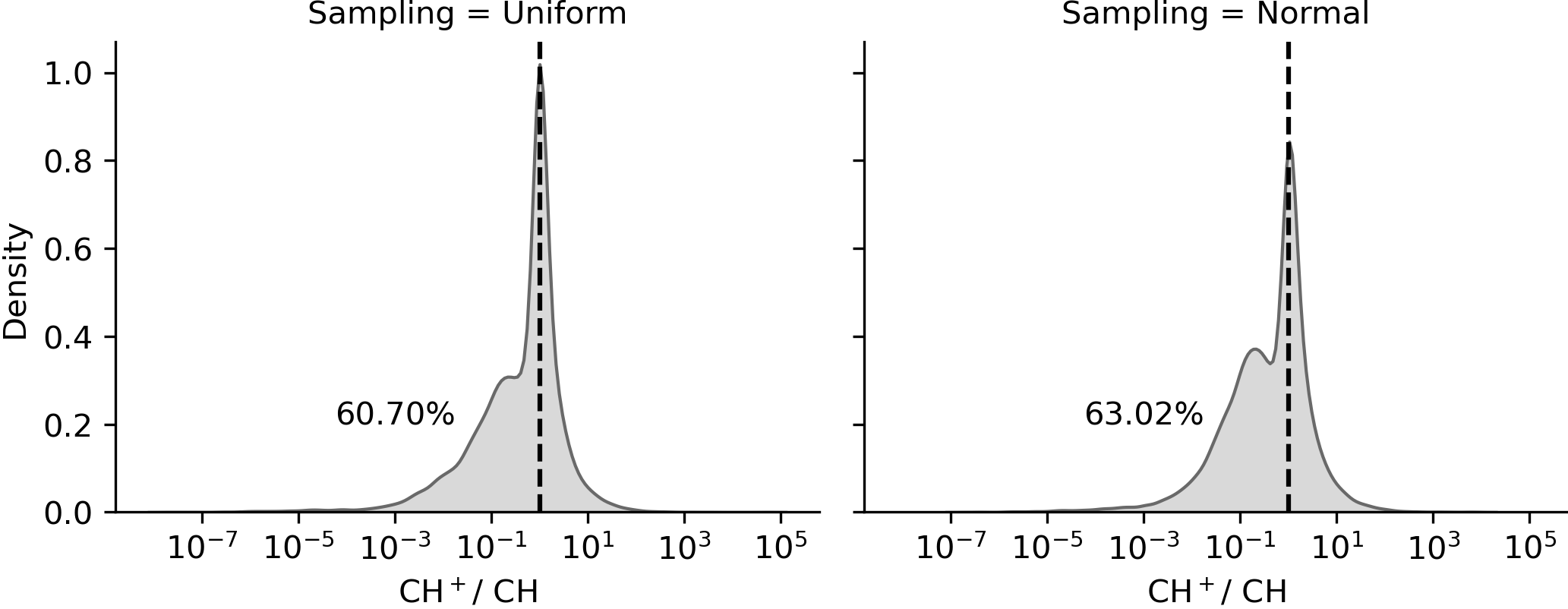}
\caption{Empirical distributions of the ratio of the function value error of $\CH^+$ divided by the function value error of $\CH$, grouped by sampling rule.}
\label{fig:funvalerr_chplus_over_ch}
\end{figure}

In Figure \ref{fig:funvalerr_chplus_over_ch}, we examine more closely
the behavior of the function value error of $\CH^+$ compared to that of
$\CH$ over all experiments. We plot the empirical distribution of the
ratio of the function value error for $\CH^+$ divided by the function
value error for $\CH$. When the ratio is less than 1.0 for a given
experiment, $\CH^+$ has a better function value error; when the ratio
is greater than 1.0, $\CH^+$ has a worse error on that experiment. In
Figure \ref{fig:funvalerr_chplus_over_ch}, the distribution of ratios is
plotted on a logarithmic scale, and two additional pieces of information
are shown. First, a vertical dotted line is plotted to mark the ratio
1.0. Second, the percentage of experiments to the left of 1.0, i.e.,
when $\CH^+$ is performing better, is annotated. We see in particular
that $\CH^+$ achieves a better function value error than $\CH$ on over
 60\% of experiments. In addition, the left tail of the distributions
show that $\CH^+$ can reduce the error by up to a factor of 1,000,
whereas the right tail indicates that the error from $\CH^+$ is usually
no more than 100 times the error from $\CH$. Finally, we note that in
both panels of Figure \ref{fig:funvalerr_chplus_over_ch}, the mode is
1.0, indicating that the most common situation is for $\CH^+$ and $\CH$
to yield the same function value error. Although we do
not show them here, analogous plots comparing the function value error
of $\CH^+$ compared to that of {\sc Box} and respectively {\sc IsoFor}
showed similar results.

 In the online supplemental material, we further provide sensitivity analysis of our method with respect to certain design choices, such as the noise level $\sigma$, the sample size $N$, and the sampling rules, e.g., {\em Normal\/} versus {\em Uniform\/}.

\subsection{Results with the enlarged convex hull} \label{sec:numresults:enlarged}

Next we recreate the results of the previous subsection except that here we test the validity domains
$
V \in \left\{ \text{{\sc Box}, CH$^+$, CH$^+_{0.05}$, CH$^+_{0.10}$} \right\},
$
where CH$^+_{0.05}$ and CH$^+_{0.10}$ refer to enlarged versions of
$\CH^+$ corresponding to $\varepsilon$ values $0.05$ and $0.10$,
respectively. The purpose of these experiments are to highlight the
performance of the enlarged convex hulls relative to $\CH^+$ and {\sc
Box}, and the results are shown in Table \ref{tab:all_enlarged}. Note
that the results for {\sc Box} are scaled to 1.00 as before, and the
rows for $\CH^+$ in Table \ref{tab:all_enlarged} match the rows for
$\CH^+$ in Table \ref{tab:all}, i.e., Table \ref{tab:all_enlarged}
contains some repeated, relevant information from Table \ref{tab:all}.

We also mention that, in actuality, the $\varepsilon$ values are $0.05$ and $0.10$ times a scaling factor that is proportional to the diameter of the ground set $X$, which is an $n_1$-dimensional box in each case. Specifically, we use the scale factor diam($X$), the diameter of $X$, to account for the natural size of the domain.

Examining the function value errors in Table \ref{tab:all_enlarged}, we
see that in nearly all groups of four errors (each group corresponding
to the four validity domains), as $\varepsilon$ increases from 0 to 0.05
to 0.10, we see a corresponding increase in the function value error. In
words, increased flexibility of the enlarged convex hull typically comes
with a sacrifice in the function value error. Examining the optimal
value and optimal solution errors, we do not see a consistent trend,
perhaps because $\CH^+$ and its enlarged version $\CH^+_{\varepsilon}$
are designed to manage the function value error and hence there is
naturally a weaker link with the optimal value and optimal solution
errors.

% latex table generated in R 4.4.1 by xtable 1.8-4 package
% Thu Feb 20 21:24:05 2025
\begin{table}[htbp]
\centering
\begin{tabular}{ll||cc|cc|cc}
Function &
\begin{tabular}[c]{@{}c@{}} Validity \\ Domain \end{tabular} &
\multicolumn{2}{c|}{%
\begin{tabular}[c]{@{}c@{}} Median Function \\ Value Error \end{tabular}%
} &
\multicolumn{2}{c|}{%
\begin{tabular}[c]{@{}c@{}} Median Optimal \\ Value Error \end{tabular}%
} &
\multicolumn{2}{c}{%
\begin{tabular}[c]{@{}c@{}} Median Optimal \\ Solution Error \end{tabular}%
} \\ 
& & {\em Uniform\/} & {\em Normal\/} & {\em Uniform\/} & {\em Normal\/} & {\em Uniform\/} & {\em Normal\/} \\
  \hline\hline
Beale & {\sc Box} & 1.00 & 1.00 & 1.00 & 1.00 & 1.00 & 1.00 \\ 
   & $\CH^+$ & \textbf{0.09} & \textbf{0.35} & \textbf{0.16} & \textbf{0.72} & 0.86 & 0.79 \\ 
   & $\CH^+_{0.05}$ & 0.83 & 0.85 & 0.87 & 0.85 & \textbf{0.68} & \textbf{0.65} \\ 
   & $\CH^+_{0.10}$ & 0.99 & 1.03 & 1.00 & 0.96 & 0.75 & 0.81 \\ \hline
  Griewank & {\sc Box} & 1.00 & 1.00 & \textbf{1.00} & 1.00 & 1.00 & 1.00 \\ 
   & $\CH^+$ & 0.49 & \textbf{0.53} & 1.18 & 0.90 & \textbf{0.89} & 1.08 \\ 
   & $\CH^+_{0.05}$ & \textbf{0.44} & 0.69 & 1.04 & \textbf{0.82} & 0.95 & \textbf{0.56} \\ 
   & $\CH^+_{0.10}$ & 0.51 & 0.83 & 1.02 & \textbf{0.82} & 0.96 & 0.58 \\ \hline
  Peaks & {\sc Box} & \textbf{1.00} & 1.00 & 1.00 & 1.00 & 1.00 & 1.00 \\ 
   & $\CH^+$ & 1.02 & \textbf{0.68} & 1.00 & 0.94 & 0.93 & 0.92 \\ 
   & $\CH^+_{0.05}$ & 1.89 & 1.18 & 1.00 & \textbf{0.90} & \textbf{0.77} & \textbf{0.54} \\ 
   & $\CH^+_{0.10}$ & 1.88 & 1.24 & 1.00 & 1.00 & 0.78 & 0.62 \\ \hline
  Powell & {\sc Box} & 1.00 & 1.00 & 1.00 & 1.00 & 1.00 & 1.00 \\ 
   & $\CH^+$ & \textbf{0.09} & \textbf{0.17} & \textbf{0.15} & \textbf{0.23} & \textbf{0.78} & 0.63 \\ 
   & $\CH^+_{0.05}$ & 0.44 & 0.50 & 0.42 & 0.36 & 0.87 & \textbf{0.59} \\ 
   & $\CH^+_{0.10}$ & 0.86 & 0.99 & 0.82 & 0.67 & 0.93 & 0.65 \\ \hline
  Qing & {\sc Box} & 1.00 & 1.00 & 1.00 & 1.00 & 1.00 & 1.00 \\ 
   & $\CH^+$ & 0.41 & 0.45 & 0.79 & 0.95 & 0.70 & 0.54 \\ 
   & $\CH^+_{0.05}$ & \textbf{0.31} & \textbf{0.41} & 0.69 & \textbf{0.92} & 0.51 & \textbf{0.36} \\ 
   & $\CH^+_{0.10}$ & 0.48 & 0.57 & \textbf{0.67} & \textbf{0.92} & \textbf{0.50} & 0.39 \\ \hline
  Quintic & {\sc Box} & 1.00 & 1.00 & 1.00 & 1.00 & 1.00 & 1.00 \\ 
   & $\CH^+$ & \textbf{0.13} & \textbf{0.06} & 1.00 & 0.25 & 0.92 & \textbf{0.74} \\ 
   & $\CH^+_{0.05}$ & 0.30 & 0.16 & \textbf{0.70} & \textbf{0.18} & \textbf{0.86} & 0.75 \\ 
   & $\CH^+_{0.10}$ & 0.60 & 0.27 & 0.78 & 0.21 & 0.93 & 0.84 \\ \hline
  Rastrigin & {\sc Box} & 1.00 & 1.00 & \textbf{1.00} & \textbf{1.00} & 1.00 & 1.00 \\ 
   & $\CH^+$ & \textbf{0.68} & \textbf{0.49} & 1.13 & 1.54 & 0.70 & \textbf{0.63} \\ 
   & $\CH^+_{0.05}$ & 0.77 & 0.71 & 1.07 & 1.31 & \textbf{0.66} & 0.77 \\ 
   & $\CH^+_{0.10}$ & 0.82 & 0.82 & 1.04 & 1.13 & 0.70 & 0.82 \\ 
\end{tabular}
\caption{Errors for all 37,800 experiments, grouped by function, type of error, and sampling rule. In each group of four errors corresponding to four validity domains, the errors are scaled so that {\sc Box} has value 1.00 and the smallest error is shown in a bold font in order to facilitate comparison among the different validity domains.}
\label{tab:all_enlarged}
\end{table}

\subsection{Results for larger dimensions}

The computational scalability of any CL approach faces several
obstacles, which may grow more severe with the base dimension $n_1$,
the sample size $N$, the complexity of the learned model ${\cal M}$,
and the resultant size and structure of the mixed-integer programming
model. Moreover, each validity domain will scale differently than other
domains.

Even though CL will require more time and memory for larger, more
complex problems, one would hope the quality and effectiveness
of the approach does not simultaneously degrade. Accordingly, in
this subsection, we test $\CH^+$ on larger problems than in Tables
\ref{tab:all} and \ref{tab:all_enlarged}. In particular, we test an
additional function for dimension $n_1 = 20$:
\[
f(x) := \sum_{i = 1}^n \left( x_1 + \cdots + x_i \right)^2,
\]
where each $x_i \in [-200, 200]$. This function is called the {\em rotated hyper ellipsoid\/} function in \cite{LONG2019108}. 

In terms of computational scalability, for these experiments, we found that whenever ${\cal R} = \text{{\em Uniform\/}}$,
$$
{\cal M} \in
\left\{
\text{\em RandomForestRegressor, GradientBoostingRegressor\/} \right\},
$$ or $V = \text{{\sc IsoFor}}$,
%$
%{\cal R} \in \left\{ \text{\em Uniform} \right\},
%$
%$
%{\cal M} \in \left\{ \text{\em RandomForestRegressor, GradientBoostingRegressor\/} \right\}
%$ or
%$
%V \in \left\{ \text{{\sc IsoFor}} \right\},
%$
the run time of each experiment was prohibitive for our testing, e.g., in such cases, a single experiment could take over 30 minutes. So in this subsection, we tested all other combinations, i.e.,
\begin{align*}
{\cal R} &\in \{ \text{\em Normal} \}, \\
N &\in \bigl\{ 500, 1000, 1500 \bigr\}, \\
s &\in \bigl\{ 2023, 2024, \dots, 2122 \bigr\} \\
\sigma &\in \bigl\{0.0, 0.1, 0.2 \bigr\} \\
{\cal M} &\in \left\{ \text{\em MLPRegressor\/} \right\} \\
V &\in \left\{ \text{{\sc Box}, CH, CH$^+$, CH$^+_{0.05}$, CH$^+_{0.10}$} \right\}.
\end{align*}
This amounted to 900 experiments with five validity domains tested
in each experiment. We also mention that, for {\em
MLPRegressor\/}, we chose a hidden layer size of 50 for these tests,
compared to 30 in the results cataloged in Tables \ref{tab:all} and
\ref{tab:all_enlarged}. All other details of our {\em MLPRegressor\/}
implementation remained the same.

The results are shown in Table \ref{tab:all_lareger_n}, which displays similar patterns as seen in the previous Tables \ref{tab:all} and \ref{tab:all_enlarged}. In other words, we see consistent performance of $\CH^+$ and $\CH^+_\varepsilon$ in larger dimensions for this function.

% latex table generated in R 4.4.1 by xtable 1.8-4 package
% Thu Feb 20 21:24:06 2025
\begin{table}[ht]
\centering
\begin{tabular}{ll||c|c|c}
Function &
\begin{tabular}[c]{@{}c@{}} Validity \\ Domain \end{tabular} &
\begin{tabular}[c]{@{}c@{}} Median Function \\ Value Error \end{tabular} &
\begin{tabular}[c]{@{}c@{}} Median Optimal \\ Value Error \end{tabular} &
\begin{tabular}[c]{@{}c@{}} Median Optimal \\ Solution Error \end{tabular} \\  
  \hline
Rotated Hyper Ellipsoid & {\sc Box} & 1.00 & 1.00 & 1.00 \\ 
   & $\CH$ & 0.09 & 0.12 & 0.26 \\ 
   & $\CH^+$ & \textbf{0.01} & \textbf{0.03} & \textbf{0.21} \\ 
   & $\CH^+_{0.05}$ & 0.05 & 0.06 & 0.25 \\ 
   & $\CH^+_{0.10}$ & 0.10 & 0.15 & 0.32 \\
\end{tabular}
\caption{Errors for all 900 experiments, grouped by function, type of error, and sampling rule. In each group of five errors corresponding to five validity domains, the errors are scaled so that {\sc Box} has value 1.00 and the smallest error is shown in a bold font in order to facilitate comparison among the different validity domains.}
\label{tab:all_lareger_n}
\end{table}

\section{Two Stylized Optimization Models} \label{sec:stylized}

In this section, we examine the performance of our extended validity domain $\CH^+$ in the context of two stylized optimization problems. In both cases, we see that $\CH^+$ and its enlarged variants more effectively manage the function value and feasibility errors.

\subsection{A simple nonlinear optimization}

Consider the true optimization problem
$
v^* := \min_{x \in \mathbb{R}^{n_1}} \left\{ c^T x : \|x\| \le 1, x \in [-1,1]^{n_1} \right\}
$
where $c \in \mathbb{R}^{n_1}$ is an arbitrary vector satisfying $\|c\| = 1$. It is easy to see that $v^* = -1$ and the unique optimal solution is $x^* = -c$. This is an instance of (\ref{equ:true}) with $f(x) := c^T x$, $X := [-1,1]^{n_1}$, $g$ nonexistent, $\theta(y) := y - 1$, and $h(x) := \|x\|$. In this subsection, we will consider small values of $n_1$, specifically $n_1 \in \{5, 10\}$.

We conduct a single experiment by randomly generating $c$ uniformly on the surface of the Euclidean unit ball, sampling $N = 1,000$ points $x$ such that the Euclidean norm of $x$ is uniform in $[0.5, 1.5]$, and evaluating $h(x) + 0.05 \epsilon$ on all samples, where $\epsilon$ is a standard normal random variable. We train the function $\hat h$ using the Python function {\em MLPRegressor\/} just as in Section \ref{sec:numresults} using the same hyperparameter choices.

For a single experiment, we then test the validity domains {\sc Box}, $\CH$, $\CH^+$, $\CH^+_{0.05}$, and $\CH^+_{0.10}$, where, as in Section \ref{sec:numresults}, the $\varepsilon$ values are in fact scaled by the diameter of the ground set $X$. We also attempted to test {\sc IsoFor}, but its overall run times were on the order of hours per experiment compared with seconds. So we decided not to test {\sc IsoFor} in order to streamline our computational testing.

$\CH^+$ is constructed as the convex hull of $F_N^+$ according to (\ref{equ:CH+}), where $F_N^+$ is the set of feasible samples in the full extended space defined by (\ref{equ:feasdataset+}). Despite the fact that the evaluations of $h$ are noisy, we include a point $x$ in $F_N^+$ if $\|x\| + 0.05 \epsilon \le 1$, where $\epsilon$ is the added noise (not to be confused with $\varepsilon$ in the definition of $\CH^+_\varepsilon$). In particular, the noise affects the construction of the set $F_N^+$. Because $x$ has been sampled with radius uniform in $[0.5, 1.5]$, the expected cardinality of $F_N^+$ is approximately $N/2$.

For $n_1 = 5$, we run 100 experiments and show the median errors in
Table \ref{tab:simple}, and then we repeat the same experiments for
$n_1 = 10$. Similar to Tables \ref{tab:all}-\ref{tab:all_lareger_n}
in Section \ref{sec:numresults}, we collect the results in Table
\ref{tab:simple}, grouped by $n_1$. For each sub-grouping of
five errors, we scale such that {\sc Box} has error
1.00. First and foremost, we see that $\CH^+$ achieves
much better function and feasibility errors than {\sc Box} and $\CH$.
Secondly, there appears to be a tradeoff between achieving good function
and feasibility errors versus achieving good optimal value and optimal
solution errors. For example, $\CH^+$ does best in terms of function
and feasibility errors but worst in terms of the optimal value and
optimal solution errors. In contrast, $\CH^+_{0.10}$ sacrifices function
and feasibility errors while improving the other two errors. In this
particular case, $\CH^+_{0.05}$ seems to strike a nice balance between
all error types.

%Results for n1 = 5
%\begin{tabular}{rrrr}
%1.000000 & 1.000000 & 1.000000 & 1.000000 \\
%0.897929 & 1.000000 & 0.938698 & 0.897929 \\
%0.175101 & 1.995991 & 1.000638 & 0.000000 \\
%0.206860 & 0.863530 & 0.860680 & 0.000000 \\
%0.300952 & 0.382543 & 0.712615 & 0.242263 \\
%\end{tabular}

%Results for n1 = 10
%\begin{tabular}{rrrr}
%1.000000 & 1.000000 & 1.000000 & 1.000000 \\
%0.191521 & 0.135313 & 0.514009 & 0.191521 \\
%0.050723 & 1.219593 & 0.714287 & 0.000000 \\
%0.073877 & 0.766101 & 0.639375 & 0.000000 \\
%0.102352 & 0.398441 & 0.546789 & 0.061661 \\
%\end{tabular}

\begin{table}[]
\centering
\begin{tabular}{ll|c|c|c|c}
    $n_1$ &
  \begin{tabular}[c]{@{}c@{}}Validity\\ Domain \end{tabular}  & 
  \begin{tabular}[c]{@{}c@{}}Median Function\\ Value Error\end{tabular} &
  \begin{tabular}[c]{@{}c@{}}Median Optimal\\ Value Error\end{tabular} &
  \begin{tabular}[c]{@{}c@{}}Median Optimal \\ Solution Error\end{tabular} &
  \begin{tabular}[c]{@{}c@{}}Median \\ Feasibility Error\end{tabular} \\ \hline\hline
5 & {\sc Box}     & 1.00 & 1.00 & 1.00 & 1.00 \\
 & CH             & 0.90 & 1.00 & 0.94 & 0.90 \\
 & $\CH^+$        & \textbf{0.18} & 2.00 & 1.00 & \textbf{0.00} \\ 
 & $\CH^+_{0.05}$ & 0.20 & 0.86 & 0.86 & \textbf{0.00} \\ 
 & $\CH^+_{0.10}$ & 0.30 & \textbf{0.38} & \textbf{0.71} & 0.24 \\ \hline
10 & Box            & 1.00 & 1.00  & 1.00   & 1.00 \\
   & CH             & 0.19 & \textbf{0.14}  & \textbf{0.51}   & 0.19 \\
   & $\CH^+$        & \textbf{0.05} & 1.22  & 0.71   & \textbf{0.00} \\
   & $\CH^+_{0.05}$ & 0.07 & 0.77  & 0.64   & \textbf{0.00} \\ 
   & $\CH^+_{0.10}$ & 0.10 & 0.40  & 0.55   & 0.06 \\
\end{tabular}
\caption{Median errors of the {\sc Box}, $\CH$, $\CH^+$, $\CH^+_{0.05}$, and $\CH^+_{0.10}$ validity domains in the stylized simple nonlinear optimization model, grouped by dimension $n_1$. The median is taken over 100 experiments, and each sub-grouping of five errors is scaled so that {\sc Box} has error $1.00$ with a bold font used to highlight the minimum error in that sub-group.}
\label{tab:simple}
\end{table}

Since $\CH^+$ has been constructed as the convex hull of $F_N^+$,
which includes only feasible samples---and since the true feasible set
is convex---it is perhaps not surprising that $\CH^+$ achieves the
smallest feasibility error. So we repeated the same tests but for the
convex hull of $D_N^+$ defined in (\ref{equ:dataset+}), i.e., both
feasible and infeasible sample points are included.
This is in fact an example of using a convex hull $\CH'$, which is
sandwiched between $\CH$ and the full $\CH^+$, as discussed in Section
\ref{sec:new_vd_and_variants}. The results showed the same trends, in
particular that $\CH'$ and its enlarged variants manage feasibility well
in this case, too.

% Results for n1 = 5
% \begin{tabular}{rrrr}
% 1.000000 & 1.000000 & 1.000000 & 1.000000 \\
% 0.897929 & 1.000000 & 0.938698 & 0.897929 \\
% 0.109684 & 0.421588 & 0.673787 & 0.058274 \\
% 0.572880 & 0.787990 & 0.749656 & 0.572880 \\
% 0.947309 & 1.000000 & 0.927507 & 0.947309 \\
% \bottomrule
% \% end{tabular}

% Results for n1 = 10
% \begin{tabular}{rrrr}
% 1.000000 & 1.000000 & 1.000000 & 1.000000 \\
% 0.191521 & 0.135313 & 0.514009 & 0.191521 \\
% 0.074842 & 0.769716 & 0.606933 & 0.000000 \\
% 0.096364 & 0.145740 & 0.450915 & 0.096364 \\
% 0.267019 & 0.250700 & 0.443078 & 0.267019 \\
% \bottomrule
% \end{tabular}

\subsection{A price optimization}

We also test $\CH^+$ on a stylized price optimization problem, which serves as a preview of the case study in Section \ref{sec:case}. Imagine a company with two substitute products, labeled as products 1 and 2. Product 1 is the low-price product, and product 2 is the high-price product. Demands $d_1$ and $d_2$ for the respective products are functions of their prices, $p_1$ and $p_2$. The company would like to set prices so as to maximize revenue under various constraints. The true model is:
\begin{subequations} \label{equ:priceopt}
\begin{align}
        \max_{p_1, p_2} \quad &d_1 p_1 + d_2 p_2 \label{equ:priceopt:obj}\\
        \st \quad \ &p_1 \in [6.5, 9.5] \quad 
                  p_2 \in [7.5, 10.5] \quad 
                  p_1 + 0.5 \le p_2 \le p_1 + 1.5 \label{equ:priceopt:price}\\
                  &d_1 = 10^7 p_1^{-3.2} p_2 ^{0.5} \quad 
                  d_2 = 10^7 p_1^{1.5} p_2 ^{-2.2} \quad 
                  d_1 + d_2 \le 1.8 \times 10^6. \label{equ:priceopt:demand}
    \end{align}
\end{subequations}
Here, the objective function (\ref{equ:priceopt:obj}) is the total revenue of the company; 
constraint (\ref{equ:priceopt:price}) describes the allowable prices; and constraint (\ref{equ:priceopt:demand}) describes the demand functions for products 1 and 2 as well as a cap on demand corresponding what can actually can be sold. Gurobi solves this problem efficiently, reporting a true optimal solution of $(p_1^*, p_2^*) = (9.50, 10.31)$. At optimality, the demand constraint is active.

To perform numerical tests, we replace the demand functions with learned
models and enforce different validity-domain constraints. To learn the
demand functions, we uniformly sampled 1,000 points in the box domain of
$(p_1, p_2)$, evaluated the corresponding $(d_1, d_2)$ without noise,
and learned the approximation $(\hat{d}_1, \hat{d}_2)$ using a simple
quadratic regression. (Note that if the data had been learned with a
log-log regression, then the learned model would be exact.) We tested
{\sc Box}, $\CH$, {\sc IsoFor}, $\CH^+$,
$\CH^+_{0.05}$, and $\CH^+_{0.10}$. The convex hull used for $\CH$ is
taken over the 2-dimensional samples $(p_1, p_2)$, while our extended
convex hull $\CH^+$ is taken over the 4-dimensional samples $(p_1,
p_2, d_1, d_2)$. In particular, we make use of both
feasible and infeasible samples, and we do not append the objective
value. As such, we actually implement a variant $\CH'$ of $\CH^+$,
which is sandwiched between $\CH$ and $\CH^+$, as discussed in Section
\ref{sec:new_vd_and_variants}. However, in our discussion here, we keep
the name $\CH^+$ for simplicity.

The median errors of the validity
domains are summarized in Table \ref{tab:price}, where as before we
scale the {\sc Box} errors to 1.00. We see clearly that $\CH^+$ achieves
much lower function value and feasibility errors than either {\sc Box}
or $\CH$, although the optimal value and optimal solution errors are
notably higher. We believe that overall this constitutes an advantage
of $\CH^+$ and its variants because, generally speaking,
one cannot have true optimality without true feasibility.
Regarding the enlarged versions of $\CH^+$, while we have chosen
$\varepsilon = 0.05$ and $\varepsilon = 0.10$ to be consistent with
the experiments in prior sections of the paper, it would appear that
both enlarged validity domains achieve the same errors as {\sc Box},
indicating that these values of $\varepsilon$ are too loose in this
particular case.

\begin{table}[]
\centering
\begin{tabular}{l|cccc}
   & 
  \begin{tabular}[c]{@{}c@{}}Median Function\\ Value Error\end{tabular} &
  \begin{tabular}[c]{@{}c@{}}Median Optimal\\ Value Error\end{tabular} &
  \begin{tabular}[c]{@{}c@{}}Median Optimal \\ Solution Error\end{tabular} &
  \begin{tabular}[c]{@{}c@{}}Median \\ Feasibility Error\end{tabular} \\ \hline
{\sc Box}      & 1.00 & \textbf{1.00}  & \textbf{1.00}   & 1.00 \\
CH             & 0.99 & 1.06  & 1.07   & 0.99 \\
{\sc IsoFor}   & 0.67 & 3.81  & 4.97  & 0.66 \\
$\CH^+$        & \textbf{0.16} & 10.81 & 16.92 & \textbf{0.11} \\
$\CH^+_{0.05}$ & 1.00 & \textbf{1.00}  & \textbf{1.00}  & 0.99 \\
$\CH^+_{0.10}$ & 1.00 & \textbf{1.00}  & \textbf{1.00}  & 1.00 \\
\end{tabular}
\caption{Median errors of all validity domains in the price optimization model. The median is taken over 10,000 experiments, and each sub-grouping of three errors is scaled so that {\sc Box} has error $1.00$ with a bold font used to highlight the minimum error in that sub-group.}
\label{tab:price}
\end{table}

\section{A Case Study} \label{sec:case}

In this section, we investigate an \href{https://github.com/Gurobi/modeling-examples/blob/master/price_optimization/price_optimization_gurobiML_wls.ipynb}{avocado-price optimization model} recently described by Gurobi Optimization, the makers of \cite{gurobi}. The goal is to set the prices and supply quantities for avocados across eight regions of the United States while incorporating transportation and other costs and maximizing the total national profit from avocado sales. A critical component of the model is the relationship between avocado prices and the demand for avocados in each region. Gurobi Optimization proposed to learn the demand function using ML techniques based on observed sales data.  Then the demand function could be embedded in the optimization model, thus creating a case study for constraint learning (CL). We revisit this study in light of our extended validity domain $\CH^+$ proposed in Section \ref{sec:extvdom}.

Note that, in this case study, there is no ground truth, and so we are unable to measure the errors defined in Section \ref{sec:errors}. Instead, we seek experimental insights from this case.

\subsection{Avocado price model}

The eight regions are indexed by $r\in\{1,\ldots,8\}$, and the total units of avocado imported into a single port in the United States is denoted by $B$. The per-unit cost of waste is $\alpha$, which is independent of the region, and the transportation cost from the port to region $r$ is $\beta_r$. The learned demand function for avocados in region $r$ at sales price $p$ is denoted $\hat d(p,r)$. The optimization variables are the unit price $p_r$ and units supplied $x_r$ for region $r \in \{1,\ldots,8\}$. Auxiliary variables $s_r$ and $w_r$ represent the units sold and units wasted per region, respectively. The model formulation is
\begin{subequations} \label{avocado}
\begin{align}
    \max \quad &\sum_{r=1}^8 \left( p_r s_r- \alpha \, w_r-\beta_r x_r \right) \label{avocado:obj} \\
    \st \quad \, &\sum_{r=1}^8 x_i=B  \label{avocado:supply} \\
    &s_r \leq \min\{x_r, \hat d(p_r,r)\} \quad w_r = x_r-s_r \quad p_r \ge 0, \ x_r \ge 0 \quad \forall \ r \in\{1,\ldots,8\}. \label{avocado:bounds}
%    &s_r \leq \min\{x_r, \hat d(p_r,r)\} \quad &\forall \ r\in\{1,\cdots,8\}, \label{avocado:sales} \\ 
%    &w_r = x_r-s_r &\forall \ r\in\{1,\cdots,8\}, \label{avocado:waste} \\
%    &p_r \ge 0, \ x_r \ge 0 &\forall \ r \in\{1,\cdots,8\}. \label{avocado:bounds}
\end{align}
\end{subequations}
The nonconvex, bilinear objective (\ref{avocado:obj}) calculates profit by subtracting the cost of shipping and waste from the revenue over all regions. Constraint (\ref{avocado:supply}) ensures that the total units supplied equals the import quantity, and constraint (\ref{avocado:bounds}): defines the number of units sold in a region to be no larger than the minimum of supply and predicted demand in the region; sets the number of units wasted to be the difference between number of units supplied and the number sold; and finally enforces nonnegativity of $p_r$ and $x_r$. In practice, there may also be upper bounds on $p_r$ and $x_r$.

\subsection{Dataset and predictive model} \label{sec:case:data}

The dataset to learn $\hat d(p,r)$ has been prepared by Gurobi Optimization from two sources, and the combined dataset is hosted at \href{https://github.com/Gurobi/modeling-examples/blob/master/price_optimization/HAB_data_2015to2022.csv}{Gurobi Optimization's GitHub account}. One source is the \href{https://hassavocadoboard.com/}{Hass Avocado Board} (HAB), from which Gurobi Optimization has retrieved data for the years 2019 to 2022, and the second source is \href{https://www.kaggle.com/datasets/timmate/avocado-prices-2020}{Kaggle}, which hosts data originating from HAB for the years 2015 to 2018. The total time span of the data is thus from 2015 to 2022. Note that, since true demand is not observed directly, sales data is used instead as the best available proxy of demand, and we do not consider other issues such as non-stationarity of demand over time.

The data has been cleaned and processed such that each observation includes a date (indicating the start of a calendar week), a seasonality indicator for that specific week (corresponding to $1$ for peak and 0 for off-peak), the region, the number of units of avocados sold (in millions of units), and the average price per unit. For this data set, the average units sold over all observations is 3.9 million, the average price is \$1.14, and the average revenue is \$4.23 million; these values describe a typical week in a typical region. Aggregating up to the entire United States, in a typical week, the average units sold is 31.0 million at an average price of \$1.14 for an average revenue of \$33.85 million. Restricting to just off-peak data, in a typical off-peak week, the average units sold across the U.S. is 28.6 million at an average price of \$1.14 for an average revenue of \$31.10 million.

Following Gurobi Optimization's example, we use a gradient boosting regressor implemented using the Python's {\em scikit-learn\/} package to learn the demand function $\hat d(p,r)$ based on four features of the data: year (not the specific week), the seasonality indicator, region, and price. For convenience, we label these features {\em Year\/}, {\em Peak\/}, {\em Region\/}, and {\em Price\/}, respectively. The label for the response variable in the data is {\em Units Sold\/}.

\begin{figure}[htbp]
\centering
\includegraphics[width=1.0\linewidth]{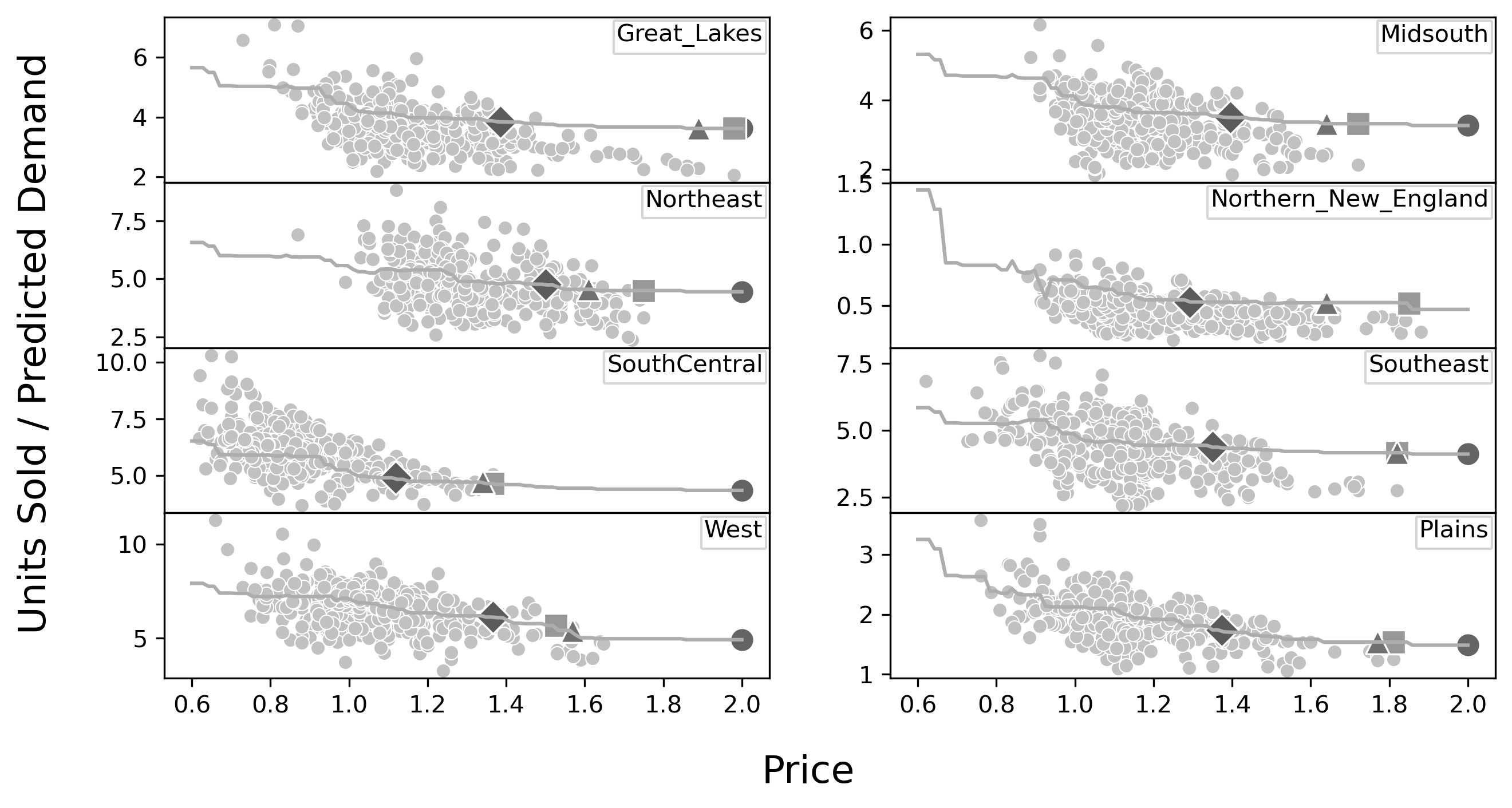}
\caption{Avocado data (light circle dots), predicted demand function (curve), and four optimal solutions. The default (i.e., no validity domain) is circle, {\sc Box} is square, $\CH$ is triangle, and $\CH^+$ is diamond.}
\label{fig:avocad_chplus}
\end{figure}

Figure \ref{fig:avocad_chplus} depicts eight panels, one for each region in the data. Here, the regions are labeled with their descriptive names, e.g., {\em Northeast\/} and {\em Plains\/}, as opposed to their index numbers $r \in \{1,\ldots,8\}$. Each panel depicts a scatter plot of light-color, circular dots, which show the observed data, where the horizontal axis is {\em Price\/} and the vertical axis is {\em Units Sold\/} (in millions). We note that each panel is plotted over the region $0.6 \le p_r \le 2.0$, representing a unit price between \$0.60 and \$2.00 for all regions. In addition, in each panel the learned function $\hat d(p_r, r)$ with fixed values {\em Year\/} $= 2023$ and {\em Peak\/} = $0$ is plotted as a curve through the data; the vertical axis is labeled as {\em Predicted Demand\/}. (Each panel also contains additional information in the form of four additional markers, which we explain in the next subsection.) We observe the expected inverse relationship between price and demand. Further, the fitted demand function does a reasonably good job capturing the relationship between price and demand in each of the eight regions.

\subsection{Effect of validity domains} \label{sec:case:vdomain}

We now solve (\ref{avocado}) for various validity domains. All experiments are executed using Gurobi 12.0 as the optimization solver. In these optimizations, we fix {\em Year\/} $= 2023$, i.e., predicting into the next year beyond the data, and we also fix {\em Peak\/} $=0$. Hence, we are optimizing prices for a single off-peak week in 2023. We set the remaining model parameters as Gurobi Optimization has done with $B = 30$, $\alpha = 0.1$, and the following transportation coefficients $\beta_r$: Great Lakes = 0.3, Midsouth = 0.1, Northeast = 0.4, Northern New England = 0.5, South Central = 0.3, Southeast = 0.2, West = 0.2, and Plains = 0.2. Gurobi was able to solve all instances to its default optimality tolerance in just a few minutes on a modern laptop.

When (\ref{avocado}) is optimized with a uniform upper bound of $p_r \le 2.0$ for all regions, the optimal value is \$45.84 million. A corresponding best solution is depicted with circle markers in the panels of Figure \ref{fig:avocad_chplus}. In particular, for all but the Northern New England region, the best solution sets $p_r = 2.0$, i.e., the price is at the upper bound of \$2.00. Particularly striking is that the solution is visually quite far from the observed data. One may ask if this solution is trustworthy given that the demand predictions are likely to be less reliable away from the observed data.

We next test the three validity domains {\sc Box}, CH, and CH$^+$. The results for {\sc IsoFor} and $\CH^+_\varepsilon$ are included in the online supplemental material. Note that all validity domains are considered with respect to the prices $p$ but not the shipped quantities $x$ because the historical data on $x$ are not available. Our hope is that these validity domains can help alleviate the uncertainty inherent in the default solution just mentioned. We would also like to compare and contrast these three validity domains in this setting. Our observations are summarized as follows:

\begin{itemize}

\item {\sc Box}: For each $r$ independently, we constrain $p_r$ to be within its observed minimum and maximum values; see (\ref{equ:box}). The best reported value is \$39.51 million, and the corresponding prices per region are shown in Figure \ref{fig:avocad_chplus} as square markers. Compared to the circle markers of the default solution, the square markers are ``pulled back'' much closer to the data. It is clear that {\sc Box} does not allow as much extrapolation in the optimal solution, and hence one can expect the predicted demand to be more accurate.  Hence, one can expect the final profit number to be more reliable. 

\item CH: This validity domain is defined to be the convex hull in $\mathbb{R}^8$ of the 8-tuples $(p^{(i)}_1, \ldots, p^{(i)}_8)$, where $i = 1, \ldots, N$ indexes over all observations; see (\ref{equ:CH}). The value is \$38.02 million, and Figure \ref{fig:avocad_chplus} shows the prices as triangle markers. As with {\sc Box}, the prices are visually closer to observed price data, and hence one can expect the overall optimization result to be more reliable.

\item $\CH^+$: Finally, we consider  our extended validity domain, which enforces that the concatenated variables and demand predictions $( p_1, \ldots, p_8, \hat d(p_1, 1), \ldots, \hat d(p_8, 8) )$ lie in the convex hull of the observed data $( p^{(i)}_1, \ldots, p^{(i)}_8, d^{(i)}_1, \ldots, d^{(i)}_8 )$ where $i = 1,\ldots,N$; see (\ref{equ:CH+}). This convex hull lies in $\mathbb{R}^{16}$. In particular, because the historical shipped quantities $x$ are not available in the historical dataset, our $\CH^+$ here is actually a variant of the form $\CH'$ discussed in Section \ref{sec:new_vd_and_variants}. The optimized prices are shown as diamond markers in Figure \ref{fig:avocad_chplus}, and the  value is \$32.54 million. Although this value is significantly less than the preceding objective values (compare to the value of \$38.02 for $\CH$, for example), the positions of the diamonds is considerably closer to the original data set. In fact, based on the two-dimensional panels in Figure \ref{fig:avocad_chplus}, it appears that the optimal prices are actually embedded inside the data, although this may be a visual artifact of the projection of a 16-dimensional image down to eight individual 2-dimensional scatter plots. In any case, one can expect that the predictions $\hat d(p,r)$ are the most reliable for this validity domain, making the overall optimization more reliable.

\end{itemize}

As a final comment, we recall that the empirical data shows an average revenue of \$31.10 million in a typical off-peak week; see the discussion in Section \ref{sec:case:data}. The final optimal value of \$32.54 million for $\CH^+$ is certainly in line with this empirical average and constitutes an increase of 4.6\%.

\section{Conclusions} \label{sec:conclusion}

We have studied the use of validity domains to reduce the errors associated with the constraint-learning (CL) framework. Based on the intuition of using the convex hull to learn both the data set and the optimization problem itself, we have proposed a new extended validity domain, called $\CH^+$. Our numerical studies have shown that $\CH^+$, compared to other common methods in the literature, is competitive in terms of computational effort and tends especially to reduce the function value and feasibility errors. Beyond stylized numerical results, the avocado case study has shown the applicability and adaptability of our approach to real-world situations.

 While we would certainly advocate the use of $\CH^+$ as a valuable option in many situations (e.g., to reduce the function value and feasibility errors, as just mentioned), it is of course natural to ask which validity domain is best suited for any specific situation. For example, if one is primarily interested in estimating the true optimal value $v^*$, then choosing the validity domain that minimizes optimal value error would be the best choice. At this time, we unfortunately lack a solid set of practical guidelines for choosing the best validity domain in common situations. Nevertheless, we hope that this research contributes to the growing literature on how to make the best use of constraint learning.

We mention some opportunities for future research. First, the model
(\ref{equ:stdform}) based on $\hat h$ may be infeasible in general,
and when validity domains such as $\CH$ or $\CH^+$ are also enforced,
infeasibility will be, in a sense, even more likely. It will be
interesting to investigate the underlying properties that make
(\ref{equ:stdform}) feasible, to understand when $\CH^+$ maintains this
feasibility, and if not, then to develop alternate validity domains that
do maintain it. Second, we have used the convex hull in an extended
space in part because $\CH$ in the original space is well-studied and
possesses good properties. Of course, there exist other techniques
for creating validity domains in the original space as described in
Section \ref{sec:cl}. One future idea to explore is whether these other
techniques, similar to $\CH^+$, can also be effective in the extended
space. Third, while current validity domains in the literature are
model-agnostic, including $\CH^+$, future research could investigate
validity domains specific to the optimization model as well as the ML
model $\hat h$ used to learn the constraints. Such validity domains
might be termed ``model aware.'' Lastly, this paper has only dealt
with regression models $\hat h$. Exactly how $\CH^+$ might work in the
case of classification is an extremely interesting question, which, if
successful, could increase the versatility of the ideas in this paper.

\section*{Acknowledgments}

The authors wish to thank Paul Grigas for discussions on the relationship between constraint learning and contextual optimization.  The authors are also in debt to the anonymous associate editor and three referees for invaluable suggestions that have greatly improved the paper.

%\ifthenelse{\boolean{ijoc}}{}{\end{onehalfspace}}

%\bibliographystyle{INFORMS-IJOC-Template-6-10-2024/informs2014}
%\bibliographystyle{informs2014}
\bibliographystyle{plainnat}
\bibliography{paper}

\begin{thebibliography}{19}
\providecommand{\natexlab}[1]{#1}
\providecommand{\url}[1]{\texttt{#1}}
\expandafter\ifx\csname urlstyle\endcsname\relax
  \providecommand{\doi}[1]{doi: #1}\else
  \providecommand{\doi}{doi: \begingroup \urlstyle{rm}\Url}\fi

\bibitem[Balestriero et~al.(2021)Balestriero, Pesenti, and LeCun]{balestriero2021learninghighdimensionamounts}
Randall Balestriero, Jerome Pesenti, and Yann LeCun.
\newblock Learning in high dimension always amounts to extrapolation, 2021.
\newblock URL \url{https://arxiv.org/abs/2110.09485}.

\bibitem[Bengio et~al.(2021)Bengio, Lodi, and Prouvost]{bengio2021machine}
Yoshua Bengio, Andrea Lodi, and Antoine Prouvost.
\newblock Machine learning for combinatorial optimization: a methodological tour d’horizon.
\newblock \emph{European Journal of Operational Research}, 290\penalty0 (2):\penalty0 405--421, 2021.

\bibitem[Bergman et~al.(2022)Bergman, Huang, Brooks, Lodi, and Raghunathan]{bergman2019janos}
David Bergman, Teng Huang, Philip Brooks, Andrea Lodi, and Arvind~U. Raghunathan.
\newblock {\tt {JANOS}}: an integrated predictive and prescriptive modeling framework.
\newblock \emph{INFORMS J. Comput.}, 34\penalty0 (2):\penalty0 807--816, 2022.
\newblock ISSN 1091-9856,1526-5528.
\newblock \doi{10.1287/ijoc.2020.1023}.
\newblock URL \url{https://doi.org/10.1287/ijoc.2020.1023}.

\bibitem[Ceccon et~al.(2022)Ceccon, Jalving, Haddad, Thebelt, Tsay, Laird, and Misener]{ceccon2022omlt}
Francesco Ceccon, Jordan Jalving, Joshua Haddad, Alexander Thebelt, Calvin Tsay, Carl~D. Laird, and Ruth Misener.
\newblock O{MLT}: optimization \& machine learning toolkit.
\newblock \emph{J. Mach. Learn. Res.}, 23:\penalty0 Paper No. [349], 8, 2022.
\newblock ISSN 1532-4435,1533-7928.

\bibitem[Courrieu(1994)]{courrieu1994three}
Pierre Courrieu.
\newblock Three algorithms for estimating the domain of validity of feedforward neural networks.
\newblock \emph{Neural Networks}, 7\penalty0 (1):\penalty0 169--174, 1994.

\bibitem[De~Filippo et~al.(2018)De~Filippo, Lombardi, and Milano]{De_Filippo_2018}
Allegra De~Filippo, Michele Lombardi, and Michela Milano.
\newblock Methods for off-line/on-line optimization under uncertainty.
\newblock In \emph{Proceedings of the Twenty-Seventh International Joint Conference on Artificial Intelligence}, IJCAI-2018. International Joint Conferences on Artificial Intelligence Organization, July 2018.
\newblock \doi{10.24963/ijcai.2018/177}.
\newblock URL \url{http://dx.doi.org/10.24963/ijcai.2018/177}.

\bibitem[Fajemisin et~al.(2023)Fajemisin, Maragno, and den Hertog]{fajemisin2023optimization}
Adejuyigbe~O Fajemisin, Donato Maragno, and Dick den Hertog.
\newblock Optimization with constraint learning: a framework and survey.
\newblock \emph{European Journal of Operational Research}, 2023.

\bibitem[{Gurobi}(2023)]{gurobi}
{Gurobi}.
\newblock {Gurobi Optimizer Reference Manual}, 2023.
\newblock URL \url{https://www.gurobi.com}.

\bibitem[Kotary et~al.(2021)Kotary, Fioretto, Van~Hentenryck, and Wilder]{osti_10337537}
James Kotary, Ferdinando Fioretto, Pascal Van~Hentenryck, and Bryan Wilder.
\newblock End-to-end constrained optimization learning: A survey.
\newblock \emph{International Joint Conference on Artificial Intelligence}, 2021.
\newblock \doi{10.24963/ijcai.2021/610}.
\newblock URL \url{https://par.nsf.gov/biblio/10337537}.

\bibitem[Liu et~al.(2008)Liu, Ting, and Zhou]{liu2008isolation}
Fei~Tony Liu, Kai~Ming Ting, and Zhi-Hua Zhou.
\newblock Isolation forest.
\newblock In \emph{2008 eighth ieee international conference on data mining}, pages 413--422. IEEE, 2008.

\bibitem[Long et~al.(2019)Long, Wu, Liang, and Xu]{LONG2019108}
Wen Long, Tiebin Wu, Ximing Liang, and Songjin Xu.
\newblock Solving high-dimensional global optimization problems using an improved sine cosine algorithm.
\newblock \emph{Expert Systems with Applications}, 123:\penalty0 108--126, 2019.
\newblock ISSN 0957-4174.
\newblock \doi{https://doi.org/10.1016/j.eswa.2018.11.032}.
\newblock URL \url{https://www.sciencedirect.com/science/article/pii/S0957417418307528}.

\bibitem[Maragno et~al.(2025)Maragno, Wiberg, Bertsimas, Birbil, den Hertog, and Fajemisin]{maragno2025mixed}
Donato Maragno, Holly Wiberg, Dimitris Bertsimas, {\c{S}}~{\.I}lker Birbil, Dick den Hertog, and Adejuyigbe~O Fajemisin.
\newblock Mixed-integer optimization with constraint learning.
\newblock \emph{Operations Research}, 73\penalty0 (2):\penalty0 1011--1028, 2025.

\bibitem[Mistry et~al.(2021)Mistry, Letsios, Krennrich, Lee, and Misener]{mistry2021mixed}
Miten Mistry, Dimitrios Letsios, Gerhard Krennrich, Robert~M Lee, and Ruth Misener.
\newblock Mixed-integer convex nonlinear optimization with gradient-boosted trees embedded.
\newblock \emph{INFORMS Journal on Computing}, 33\penalty0 (3):\penalty0 1103--1119, 2021.

\bibitem[Sadana et~al.(2025)Sadana, Chenreddy, Delage, Forel, Frejinger, and Vidal]{SADANA2025271}
Utsav Sadana, Abhilash Chenreddy, Erick Delage, Alexandre Forel, Emma Frejinger, and Thibaut Vidal.
\newblock A survey of contextual optimization methods for decision-making under uncertainty.
\newblock \emph{European Journal of Operational Research}, 320\penalty0 (2):\penalty0 271--289, 2025.
\newblock ISSN 0377-2217.
\newblock \doi{https://doi.org/10.1016/j.ejor.2024.03.020}.
\newblock URL \url{https://www.sciencedirect.com/science/article/pii/S0377221724002200}.

\bibitem[Schweidtmann et~al.(2022)Schweidtmann, Weber, Wende, Netze, and Mitsos]{schweidtmann2022obey}
Artur~M Schweidtmann, Jana~M Weber, Christian Wende, Linus Netze, and Alexander Mitsos.
\newblock Obey validity limits of data-driven models through topological data analysis and one-class classification.
\newblock \emph{Optimization and engineering}, 23\penalty0 (2):\penalty0 855--876, 2022.

\bibitem[Shi et~al.(2024)Shi, Emadikhiav, Lozano, and Bergman]{doi:10.1287/ijoc.2022.0312}
Chenbo Shi, Mohsen Emadikhiav, Leonardo Lozano, and David Bergman.
\newblock Constraint learning to define trust regions in optimization over pre-trained predictive models.
\newblock \emph{INFORMS Journal on Computing}, 36\penalty0 (6):\penalty0 1382--1399, 2024.
\newblock \doi{10.1287/ijoc.2022.0312}.
\newblock URL \url{https://doi.org/10.1287/ijoc.2022.0312}.

\bibitem[Surjanovic and Bingham(2023)]{testfunctions}
Sonja Surjanovic and Derek Bingham.
\newblock Virtual library of simulation experiments: Test functions and datasets.
\newblock \url{https://www.sfu.ca/~ssurjano/optimization.html}, 2023.
\newblock Accessed: 2024-06-11.

\bibitem[Tang and Khalil(2024)]{tang2024pyepo}
Bo~Tang and Elias~B Khalil.
\newblock Pyepo: A pytorch-based end-to-end predict-then-optimize library for linear and integer programming.
\newblock \emph{Mathematical Programming Computation}, 16\penalty0 (3):\penalty0 297--335, 2024.

\bibitem[Tjeng et~al.(2019)Tjeng, Xiao, and Tedrake]{tjeng2018evaluating}
Vincent Tjeng, Kai~Y. Xiao, and Russ Tedrake.
\newblock Evaluating robustness of neural networks with mixed integer programming.
\newblock In \emph{International Conference on Learning Representations}, 2019.
\newblock URL \url{https://openreview.net/forum?id=HyGIdiRqtm}.

\end{thebibliography}

\clearpage  % Start the supplement on a new page
\appendix   % Optionally marks beginning of appendix material

% Reset counters so supplement starts at S1
\setcounter{section}{0}
\setcounter{equation}{0}
\setcounter{figure}{0}
\setcounter{table}{0}

% Use S-prefixed numbering for supplement
\renewcommand{\thesection}{S\arabic{section}}
\renewcommand{\theequation}{S\arabic{equation}}
\renewcommand{\thefigure}{S\arabic{figure}}
\renewcommand{\thetable}{S\arabic{table}}

% Title for Supplementary Material
\begin{center}
    \Large \bfseries Supplementary Material\\[1ex]
\end{center}

\section{Optimization Run Times for Sections \ref{sec:numresults:results} and \ref{sec:numresults:enlarged}}

For the experiments described in Sections \ref{sec:numresults:results}
and \ref{sec:numresults:enlarged}, we show in Tables
\ref{tab:opt_setup_times2} and \ref{tab:opt_solve_times2} the median
setup and solve times for all combinations of machine-learning
techniques and validity domains. By {\em setup time\/}, we mean the
time required to build and pass the optimization to Gurobi, including
the critical constraint which sets $y$ equal to the output of the
learned function $\hat h(x)$, where $\hat h$ has already been trained
and stored in memory. Although computation times are not the main focus
in this paper, we see that {\sc IsoFor} requires more time in general,
while $\CH^+$ can take between 1-4 times as long as {\sc Box} and
$\CH$. Furthermore, for solve times specifically, the enlarged versions
$\CH^+_\varepsilon$ take times that are more similar to {\sc Box} and
$\CH$ than to $\CH^+$ itself. While we are not exactly sure why this
occurs, this could represent a potential advantage of $\CH^+_\varepsilon$
in practice compared to the default $\CH^+$.

% latex table generated in R 4.4.1 by xtable 1.8-4 package
% Sat Apr 19 17:34:41 2025
\begin{table}[ht]
\centering
\begin{tabular}{l|rrrrrr}
ML Technique & {\sc Box} & $\CH$ & {\sc IsoFor} & $\CH^+$ & $\CH^+_{0.05}$ & $\CH^+_{0.10}$ \\ 
  \hline
{\em RandomForestRegressor} & 1.23 & 1.23 & 41.25 & 1.23 & 0.99 & 0.97 \\ 
  {\em GradientBoostingRegressor} & 1.28 & 1.27 & 41.47 & 1.27 & 1.01 & 1.01 \\ 
  {\em MLPRegressor} & 0.01 & 0.02 & 40.47 & 0.02 & 0.02 & 0.02 \\ 
\end{tabular}
\caption{Median optimization setup times (in seconds) over all 226,800 experiments, grouped by ML technique and validity domain.}
\label{tab:opt_setup_times2}
\end{table}

%\begin{table}[htbp]
%\centering%\small
%\begin{tabular}{l|rrrrrr}
%& {\sc Box} & $\CH$ & {\sc IsoFor} & $\CH^+_{0.05}$ & $\CH^+_{0.10}$ & $\CH^+$ \\ 
%  \hline
%{\em RandomForestRegressor} & 1.23 & 1.23 & 41.25 & 1.23 & 0.99 & 0.97  \\
%  {\em GradientBoostingRegressor} & 1.28 & 1.27 & 41.47 & 1.27 & 1.01 & 1.01  \\
%  {\em MLPRegressor} & 0.01 & 0.02 &40.47 & 0.02 & 0.02 & 0.02  \\
%\end{tabular}
%\caption{Median optimization setup times (in seconds) over all 226,800 experiments, grouped by ML technique and validity domain.}
%\label{tab:opt_setup_times2}
%\end{table}

% latex table generated in R 4.4.1 by xtable 1.8-4 package
% Sat Apr 19 17:38:57 2025
\begin{table}[ht]
\centering
\begin{tabular}{l|rrrrrr}
ML Technique & {\sc Box} & $\CH$ & {\sc IsoFor} & $\CH^+$ & $\CH^+_{0.05}$ & $\CH^+_{0.10}$ \\ 
  \hline
{\em RandomForestRegressor} & 0.86 & 1.31 & 3.56 & 1.48 & 0.81 & 0.81 \\ 
  {\em GradientBoostingRegressor} & 5.40 & 6.70 & 14.35 & 14.76 & 5.94 & 4.90 \\ 
  {\em MLPRegressor} & 0.09 & 0.23 & 24.24 & 0.39 & 0.45 & 0.46 \\ 
\end{tabular}
\caption{Median optimization solve times (in seconds) over all 226,800 experiments, grouped by ML technique and validity domain.}
\label{tab:opt_solve_times2}
\end{table}

%\begin{table}[htbp]
%\centering%\small
%\begin{tabular}{l|rrrrrr}
%& {\sc Box} & $\CH$ & {\sc IsoFor} & $\CH^+_{0.05}$ & $\CH^+_{0.10}$ & $\CH^+$ \\
%  \hline
%{\em RandomForestRegressor} & 0.86 & 1.31 & 3.56 & 1.48 & 0.81 & 0.81  \\
%  {\em GradientBoostingRegressor} & 5.40 & 6.70 & 14.35 & 14.76 & 5.94 & 4.90  \\
%  {\em MLPRegressor} & 0.09 & 0.23 & 24.24 & 0.39 & 0.45 & 0.46 \\
%\end{tabular}
%\caption{Median optimization solve times (in seconds) over all 226,800 experiments, grouped by ML technique %and validity domain.}
%\label{tab:opt_solve_times2}
%\end{table}

\section{Mean Errors for Section \ref{sec:numresults:results}}

Recall that Table \ref{tab:all} in Section \ref{sec:numresults:results} reports {\em median\/} error values. To complement Table \ref{tab:all}, we now present the exact same analysis except that we calculate {\em mean\/} values in Table \ref{tab:mean}. We find that the mean values reveal a similar pattern as the median values. In particular, $\CH^+$ produces very competitive results in terms of the function value error. In fact, in all but one case, $\CH^+$ has the smallest mean function value error. For the mean optimal value errors and mean optimal solution errors, $\CH^+$ is also competitive with other validity domain techniques perform slightly better in certain cases.

\begin{table}[htbp]
\centering
\begin{tabular}{ll||cc|cc|cc}
Function &
\begin{tabular}[c]{@{}c@{}} Validity \\ Domain \end{tabular} &
\multicolumn{2}{c|}{%
\begin{tabular}[c]{@{}c@{}} Mean Function \\ Value Error \end{tabular}%
} &
\multicolumn{2}{c|}{%
\begin{tabular}[c]{@{}c@{}} Mean Optimal \\ Value Error \end{tabular}%
} &
\multicolumn{2}{c}{%
\begin{tabular}[c]{@{}c@{}} Mean Optimal \\ Solution Error \end{tabular}%
} \\ 
& & {\em Uniform\/} & {\em Normal\/} & {\em Uniform\/} & {\em Normal\/} & {\em Uniform\/} & {\em Normal\/} \\
  \hline\hline
\multirow[t]{6}{*}{Beale} & {\sc Box} & 1.00 & 1.00 & 1.00 & 1.00 & 1.00 & 1.00 \\
 & $\CH$ & 0.99 & 0.80 & 1.01 & 1.09 & 1.00 & 0.86 \\
 & {\sc IsoFor} & 0.73 & 0.50 & \textbf{0.78} & \textbf{0.49} & 0.88 & 0.64 \\
 & $\CH^+$ & \textbf{0.21} & \textbf{0.31} & 0.81 & 0.77 & \textbf{0.42} & \textbf{0.52} \\
\cline{1-8}
\multirow[t]{6}{*}{Griewank} & {\sc Box} & 1.00 & 1.00 & 1.00 & 1.00 & \textbf{1.00} & 1.00 \\
 & $\CH$ & 0.57 & 0.80 & 0.82 & 0.92 & 1.07 & 1.00 \\
 & {\sc IsoFor} & 0.63 & 0.76 & \textbf{0.30} & \textbf{0.54} & 1.68 & 1.00 \\
 & $\CH^+$ & \textbf{0.38} & \textbf{0.31} & 0.71 & 0.71 & 1.25 & \textbf{0.97} \\
\cline{1-8}
\multirow[t]{6}{*}{Peaks} & {\sc Box} & 1.00 & 1.00 & 1.00 & 1.00 & 1.00 & 1.00 \\
 & $\CH$ & \textbf{1.00} & 1.00 & 0.97 & 1.00 & \textbf{1.00} & 1.00 \\
 & {\sc IsoFor} & 1.21 & 2.13 & 0.92 & \textbf{0.50} & 1.10 & 1.00 \\
 & $\CH^+$ & 1.01 & \textbf{0.59} & \textbf{0.93} & 0.85 & 1.01 & \textbf{0.99} \\
\cline{1-8}
\multirow[t]{6}{*}{Powell} & {\sc Box} & 1.00 & 1.00 & 1.00 & 1.00 & 1.00 & 1.00 \\
 & $\CH$ & 0.99 & 0.86 & 0.95 & 0.88 & 0.99 & 0.91 \\
 & {\sc IsoFor} & 0.90 & 0.63 & \textbf{0.69} & \textbf{0.55} & 0.99 & 0.88 \\
 & $\CH^+$ & \textbf{0.50} & \textbf{0.29} & 0.79 & 0.58 & \textbf{0.63} & \textbf{0.66} \\
\cline{1-8}
\multirow[t]{6}{*}{Qing} & {\sc Box} & 1.00 & 1.00 & 1.00 & 1.00 & 1.00 & 1.00 \\
 & $\CH$ & 0.87 & 0.44 & 0.80 & 0.75 & 1.00 & 1.00 \\
 & {\sc IsoFor} & 1.23 & 0.49 & \textbf{0.66} & 0.58 & 1.00 & 1.00 \\
 & $\CH^+$ & \textbf{0.64} & \textbf{0.28} & 0.68 & \textbf{0.56} & \textbf{0.91} & \textbf{0.96} \\
\cline{1-8}
\multirow[t]{6}{*}{Quintic} & {\sc Box} & 1.00 & 1.00 & 1.00 & 1.00 & 1.00 & 1.00 \\
 & $\CH$ & 0.86 & 0.54 & 0.89 & 0.65 & \textbf{0.95} & 0.66 \\
 & {\sc IsoFor} & 0.33 & 0.23 & \textbf{0.69} & \textbf{0.48} & \textbf{0.95} & \textbf{0.59} \\
 & $\CH^+$ & \textbf{0.20} & \textbf{0.09} & 0.79 & 0.49 & 1.00 & 0.60 \\
\cline{1-8}
\multirow[t]{6}{*}{Rastrigin} & {\sc Box} & 1.00 & 1.00 & 1.00 & 1.00 & \textbf{1.00} & \textbf{1.00} \\
 & $\CH$ & 0.79 & 0.37 & 0.84 & 0.53 & 1.02 & 1.10 \\
 & {\sc IsoFor} & 0.84 & 0.45 & 0.82 & 0.52 & 1.03 & 1.08 \\
 & $\CH^+$ & \textbf{0.69} & \textbf{0.33} & \textbf{0.71} & \textbf{0.50} & 1.09 & 1.18 \\
\cline{1-8}
\end{tabular}
\caption{Mean errors for all 37,800 experiments, grouped by function, type of error, and sampling rule. In each group of four errors corresponding to the four validity domains, the errors are scaled so that {\sc Box} has value 1.00 in order to facilitate comparison among the different validity domains.}
\label{tab:mean}
\end{table}

\section{Sensitivity Analysis for the Function Value Error in Sections \ref{sec:numresults:results}
and \ref{sec:numresults:enlarged}}

In this section, we investigate the sensitivity of the various validity domains studied in this paper to changes in numerical design choices, e.g., the choice of sampling rule. We are specifically concerned with the sensitivity of the function value error in the experiments of Sections \ref{sec:numresults:results}--\ref{sec:numresults:enlarged}.

\subsection{Sensitivity with respect to sampling rule}

\begin{figure}[htbp]
\centering
\includegraphics[width=0.4\linewidth]{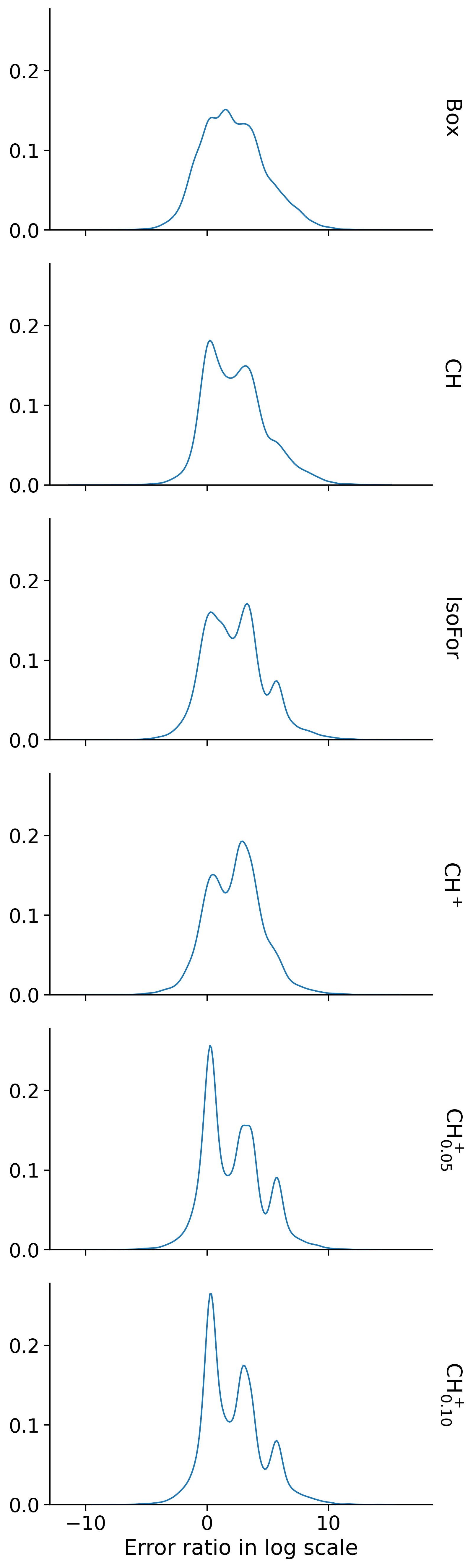}
\caption{The distribution of the ratio of the function value error obtained by Uniform divided by Normal, grouped by validity domains.}
\label{fig:sample_ratio}
\end{figure}

To investigate the sensitivity with respect to the sampling rule, i.e., {\em Uniform\/} versus {\em Normal\/}, for each 6-tuple instance of ground truth, sample size, noise level, seed, ML technique, and validity domain, we divide the function value error of {\em Uniform\/} by that of {\em Normal\/}. This ratio is greater than 1 if {\em Normal\/} exhibits better error on this instance and less than 1 if {\em Uniform\/} does. Then grouping these results by validity domain, Figure \ref{fig:sample_ratio} shows the distributions of these ratios for each validity domain. Perhaps unsurprisingly, since {\em Uniform} samples uniformly in a larger box domain, while {\em Normal} samples more closely around the global minimum in a smaller ball, we find that the ratios are skewed heavily above 1 for all validity domains, indicating that {\em Normal\/} typically exhibits less error.

\subsection{Sensitivity with respect to noise level}

Analogous to the previous subsection, we next investigate the
sensitivity of the function value error with respect to the noise level
$\sigma$, broken down by validity domain. Recall that we tested three
different values of $\sigma$, namely $0.0$, $0.1$, and $0.2$, and so for
each instance, we calculate two ratios: (i) the function value error of
$\sigma = 0.1$ divided by that of $\sigma = 0.0$; and (ii) the function
value error of $\sigma = 0.2$ divided by that of $\sigma = 0.0$. In this
sense, the noise-free case is our baseline. A ratio greater than 1 for
an instance indicates that the increased noise yields a higher function
value error compared to the no-noise case.

\begin{figure}[htbp]
\centering
\includegraphics[width=0.8\linewidth]{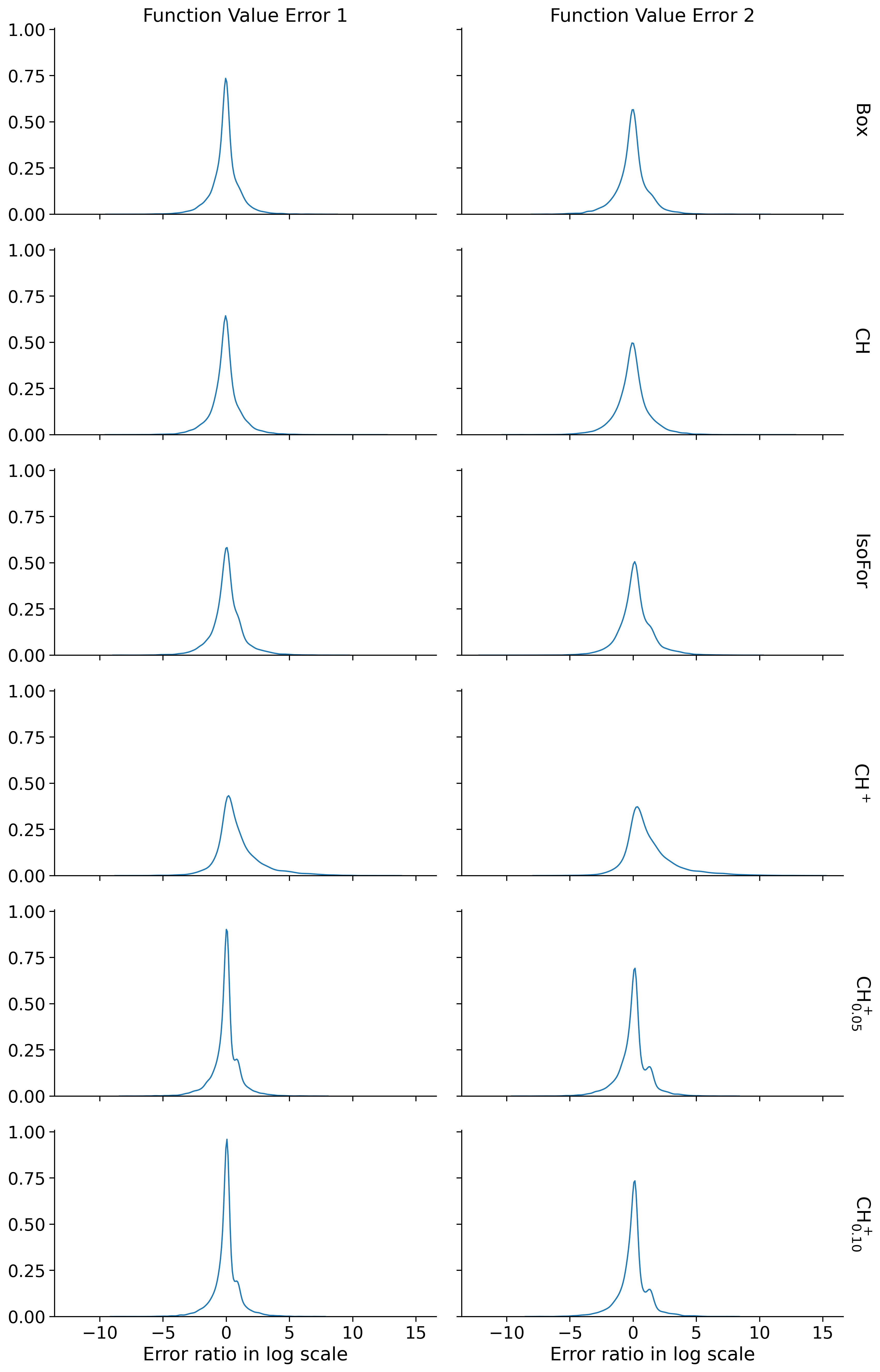}
\caption{The distribution of the ratio of the function value errors obtained by different noise levels, grouped by validity domains. Error 1 is the ratio between noise level$ = 0.1$ and noise level $=0$, Error 2 is the ratio between noise level $=0.2$, and noise level $=0$.}
\label{fig:noise_ratio}
\end{figure}

Figure \ref{fig:noise_ratio} shows the corresponding densities for each validity domain. We can see that, in general, each validity domain shows little to no sensitivity to the the noise level. This is due to the fact that each density is tightly centered around 0 with small skewnesses and large kurtosis. Table \ref{tab:noise_ratio_skew_kurt} shows the corresponding summary statistics. There is perhaps one exception, $\CH^+$, which has a larger positive skew, indicating that $\CH^+$ is more susceptible to noise, i.e., a higher noise level may slightly increase the function value errors.

%is a Facetgrid of the densities of distributions defined as follows. We first grouped the results by sampling method, validity domains, we fixed all the other parameters except the noise levels in the experiments. Then for the 3 different function value errors corresponding to 3 different noise levels, we take the ratio of the errors and plot the density of those errors in log scale. Each row corresponds to a validity domain, and each column corresponds to a ratio. The first column is the ratio between errors when noise level is 0.1 and errors when noise free, and the second column is the ratio between errors when noise level is 0,2 and errors when noise free. This is to show a more granular analysis of the experiment results with respect to the influence of different noise levels. And  shows the sknewnesses and kurtosis of the densities in figure \ref{fig:noise_ratio}, but not grouped by sampling methods. From the figure and the table, we can see that the noise level makes small differences as all the graphs are tightly centered around 0, quantitatively represented by small skewnesses and large kurtosis. There is one exception, $\CH^+$ has a more positive skewness, indicating that $\CH^+$ is more susceptible to sample noises, i.e., a higher noise level may slightly increase the errors.

\begin{table}[htbp]
\centering%\small
\begin{tabular}{l|cc|ccc}
 & \multicolumn{2}{r|}{Func Val Err 1} & \multicolumn{2}{r}{Func Val Err 2} \\
 & Skew & Kurt & Skew & Kurt \\
 \hline
Box & -0.02 & 5.36 & 0.01 & 4.07 \\
$\CH$ & 0.09 & 5.12 & 0.11 & 3.96 \\
IsoFor & 0.27 & 4.63 & 0.31 & 3.97 \\
$\CH^+$ & 1.36 & 4.34 & 1.58 & 4.62 \\
$\CH^+_{0.05}$ & -0.26 & 6.08 & -0.29 & 4.56 \\
$\CH^+_{0.10}$ & -0.16 & 6.58 & -0.21 & 4.76 \\
\end{tabular}
\caption{The skewness and kurtosis of distribution of the ratio of the errors obtained by different noise levels, grouped by validity domains. Error 1 is the ratio between noise level$ = 0.1$ and noise level $=0$, Error 2 is the ratio between noise level $=0.2$ and noise level $=0$.}
\label{tab:noise_ratio_skew_kurt}
\end{table}

\subsection{Sensitivity with respect to sample size}

\begin{figure}[htbp]
\centering
\includegraphics[width=0.8\linewidth]{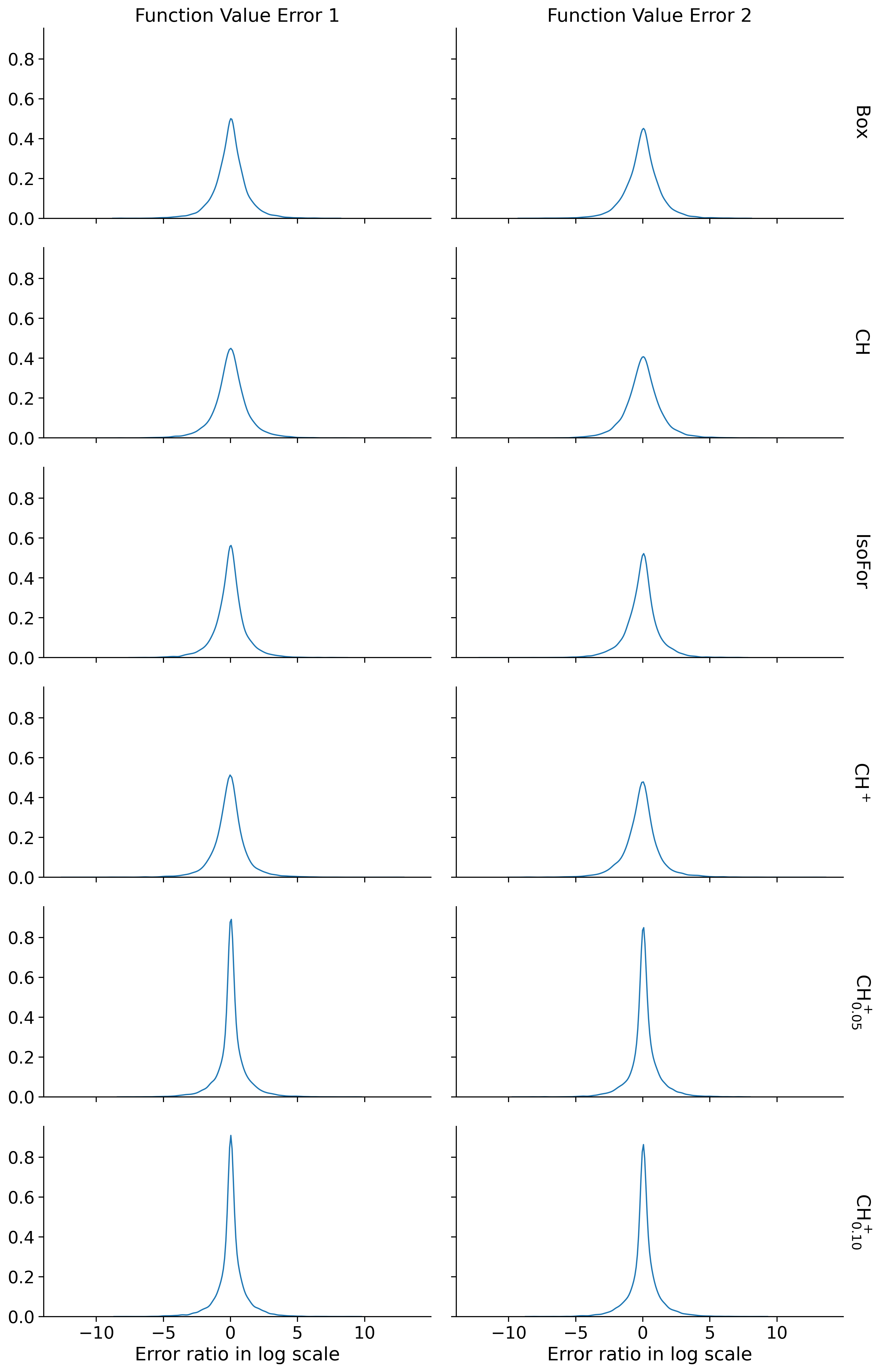}
\caption{The distribution of the ratio of the function value errors obtained by different sample sizes, grouped by validity domains. Error 1 is the ratio between sample size$ = 1000$ and sample sizes $=500$, Error 2 is the ratio between sample sizes $=1000$, and sample size $=500$.}
\label{fig:sample_sz_ratio}
\end{figure}

We next recreate the same analysis from the previous subsection but this time for the three samples sizes $N = 500, 1000, 1500$. In particular, we calculate ratios of the function value errors for $N=1000$ relative to $N=500$ and for $N=1500$ relative to $N=500$. Figure \ref{fig:sample_sz_ratio} and Table \ref{tab:sample_sz_ratio_skew_kurt} show the same patterns for $N$ as just discussed for $\sigma$. In particular, we see that changes in the sample sizes make little difference in the function value errors. Moreover, in this case, the skewness of $\CH^+$ more closely matches the skewness of the other validity domains.

\begin{table}[htbp]
\centering%\small
\begin{tabular}{l|cc|ccc}
 & \multicolumn{2}{r|}{Func Val Err 1} & \multicolumn{2}{r}{Func Val Err 2} \\
 & Skew & Kurt & Skew & Kurt \\
 \hline
Box & -0.10 & 3.60 & -0.11 & 3.24 \\
$\CH$ & 0.03 & 4.49 & 0.11 & 3.11 \\
IsoFor & 0.04 & 4.25 & -0.10 & 4.23 \\
$\CH^+$ & -0.08 & 7.17 & 0.18 & 5.89 \\
$\CH^+_{0.05}$ & 0.09 & 6.34 & 0.00 & 6.58 \\
$\CH^+_{0.10}$ & 0.05 & 6.57 & 0.02 & 6.04 \\
\end{tabular}
\caption{The skewness and kurtosis of distribution of the ratio of the errors obtained by different sample sizes, grouped by validity domains. Error 1 is the ratio between sample size$ = 1000$ and sample sizes $=500$, Error 2 is the ratio between sample sizes $=1500$, and sample size $=500$.}
\label{tab:sample_sz_ratio_skew_kurt}
\end{table}

\subsection{Sensitivity with respect to function type}

\begin{figure}[htbp]
\centering
\includegraphics[width=0.8\linewidth]{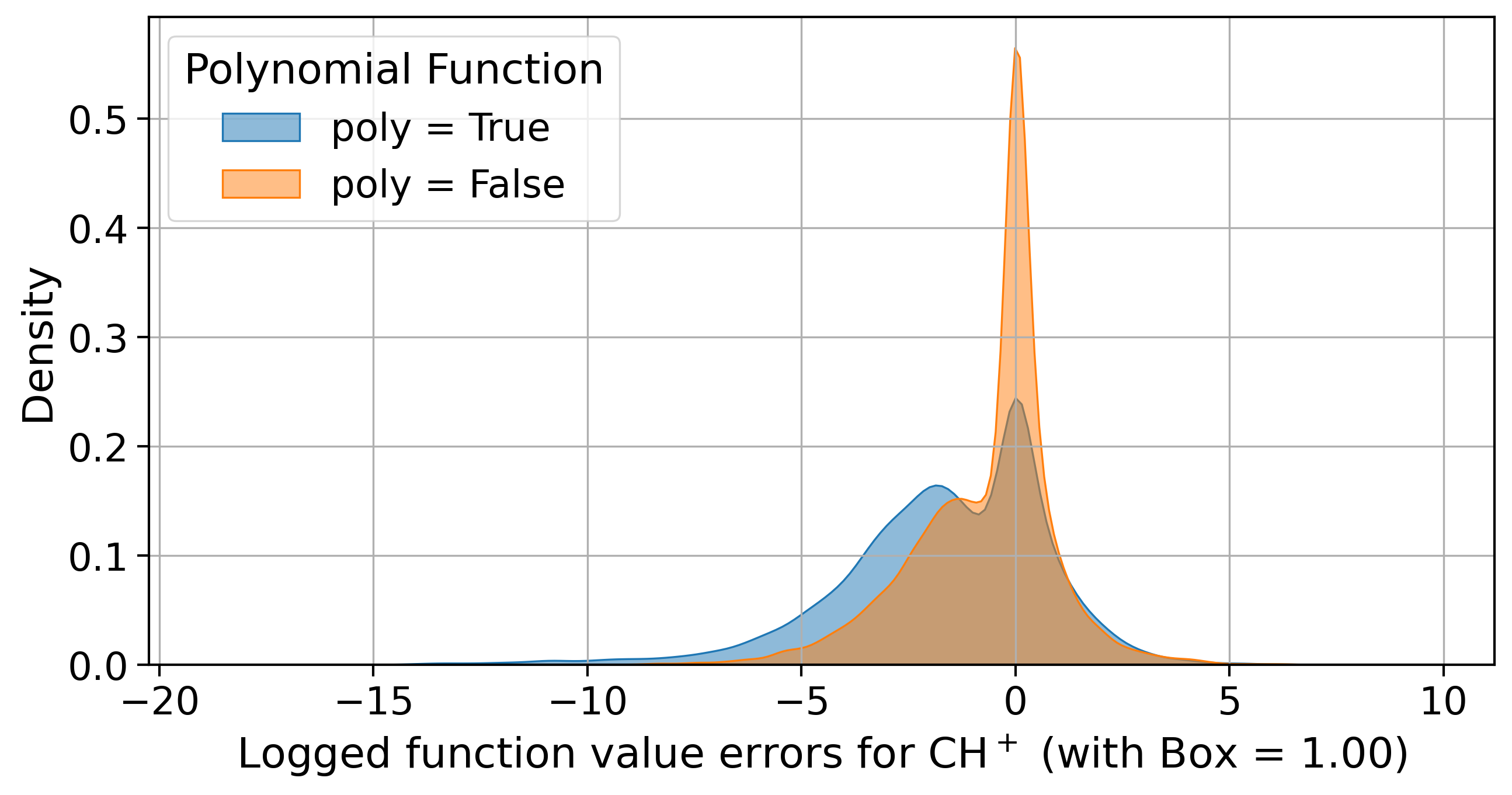}
\caption{The distribution of function value errors grouped by whether the function is a polynomial. The blue one is for polynomial functions ({\em Beale\/}, {\em Qing\/}, {\em Quintic\/}, {\em Powell\/}), and the orange one is for non-polynomial functions ({\em Peaks\/}, {\em Griewank\/}, {\em Rastrigin\/}).}
\label{fig:if_poly}
\end{figure}

Our final sensitivity analysis is with respect to whether the ground-truth function is a polynomial ({\em Beale\/}, {\em Qing\/}, {\em Quintic\/}, {\em Powell\/}) or not ({\em Peaks\/}, {\em Griewank\/}, {\em Rastrigin\/}). We note in particular that the three non-polynomial functions are based on trigonometric and exponential functions. Figure \ref{fig:if_poly} shows the densities of all function value errors in the experiments of Sections \ref{sec:numresults:results}--\ref{sec:numresults:enlarged} grouped by polynomial versus non-polynomial. We see that the distribution of the errors in the polynomial case are slightly shifted to the left compared to the non-polynomial case, which indicates that, in our experiments, polynomial functions tended to exhibit smaller function value errors.

\section{Additional Experiments for the Avocado Case Study}

\begin{figure}[htbp]
\centering
\includegraphics[width=1.0\linewidth]{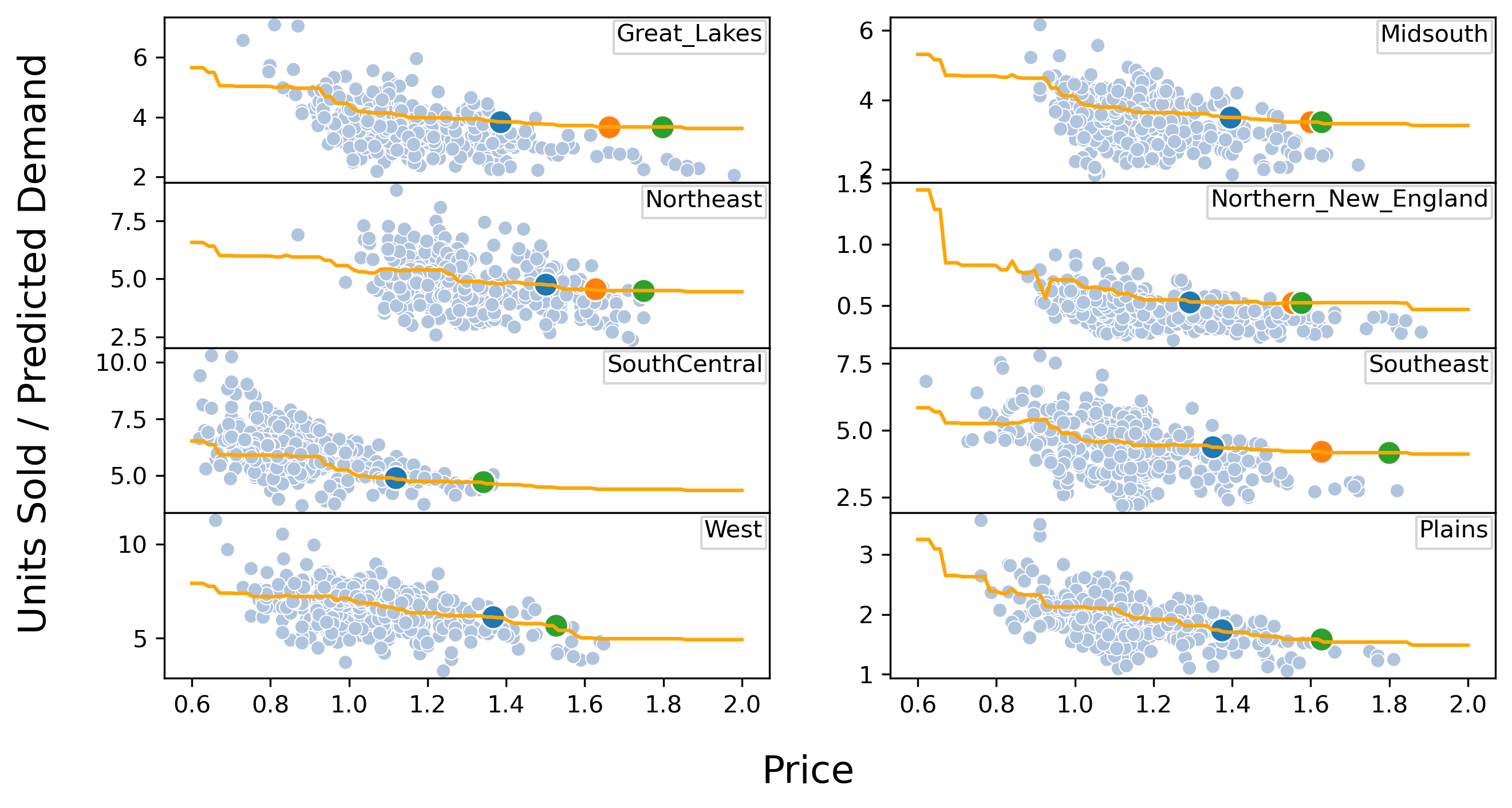}
\caption{Avocado data (light blue dots), predicted demand function (orange curve), and three optimal solutions obtained from convex hull related validity domains. The $\CH^+$ is dark blue, $\CH^+_{0.05}$ is orange, and $\CH^+_{0.10}$ is green.}
\label{fig:avocad_enlarged}
\end{figure}

Figure \ref{fig:avocad_enlarged} shows the learned optimal solutions in avocado case study when the validity domains are $\CH^+$, $\CH^+_{0.05}$ and $\CH^+_{0.10}$, where in fact the 0.05 and 0.10 values are scaled by the diameter of the box domain for this case. The learned optimal solution from $\CH^+$ is dark blue, from $\CH^+_{0.05}$ is orange, and from $\CH^+_{0.10}$ is green. We can see the trend that the smaller $\varepsilon$ is, the closer the learned optimal solution is to the center of the data points cluster. Note that in some of the graphs, the orange and green point overlap. This result show that, $\varepsilon$ is a parameter that controls the conservatism of the validity domain. We also mention that the respective optimal profits were \$32.54 million for $\CH$, \$36.54 million for $\CH^+_{0.05}$, and \$38.26 million for $\CH^+_{0.10}$.

Finally, we also ran the same experiment with the {\sc IsoFor} validity domain. It gave an optimal value of \$34.37 million, which is between $\CH^+$ and $\CH^+_{0.05}$.

%\end{document}

\end{document}